\documentclass[%
11pt,
onecolumn,
tightenlines,
superscriptaddress,
preprintnumbers,
nofootinbib,
physrev,
eqsecnum,
]{revtex4-2}

\usepackage{isomath}
\usepackage{amsmath,amsthm}
\usepackage{amsbsy}
\usepackage{amssymb}
\usepackage{amscd}
\usepackage{amsfonts}
\usepackage{stmaryrd}
\usepackage{siunitx}
\usepackage{euscript}
\usepackage[utf8]{inputenc}
\usepackage[T1]{fontenc}
\usepackage{newtxtext} 
\everymath{\displaystyle}
\usepackage{exscale}

\usepackage{graphicx}
\usepackage{boxedminipage}
\usepackage{calc}
\usepackage[usenames,dvipsnames]{xcolor}
\graphicspath{ {media/} }
\usepackage[caption=false,justification=centerlast]{subfig}

\usepackage{setspace}
\usepackage{enumitem}
\setitemize{noitemsep,topsep=0pt,parsep=0pt,partopsep=0pt}
\setenumerate{noitemsep,topsep=0pt,parsep=0pt,partopsep=0pt}
\setdescription{noitemsep,topsep=0pt,parsep=0pt,partopsep=0pt}

\usepackage{soul} 

\usepackage[,colorlinks,urlcolor=blue,citecolor=black]{hyperref}
\usepackage[normalem]{ulem}

\usepackage[small]{titlesec}

\titlespacing*{\section}{0pt}{12pt plus 4pt minus 2pt}{2pt plus 2pt minus 2pt}
\titlespacing*{\subsection}{0pt}{12pt plus 4pt minus 2pt}{2pt plus 2pt minus 2pt}
\titlespacing*\subsubsection{0pt}{12pt plus 4pt minus 2pt}{2pt plus 2pt minus 2pt}
\titlespacing*\paragraph{0pt}{12pt plus 4pt minus 2pt}{2pt plus 2pt minus 2pt}

\makeatletter

    \renewcommand*{\p@subsection}{}
    
    \renewcommand*{\p@subsubsection}{}
\makeatother

\DeclareMathOperator{\relu}{ReLU}
\DeclareMathOperator{\sigmoid}{Sigmoid}
\DeclareMathOperator{\MSE}{MSE}

\usepackage{isomath}
\usepackage{amsmath}
\usepackage{amssymb}
\usepackage{amscd}
\usepackage{amsfonts}

\newcommand{\bfx}{{\mathbold x}}

\newcommand{\bfC}{{\mathbold C}}

\newcommand{\bfR}{{\mathbold R}}

\usepackage{float}
\newcommand{\specialcell}[2][c]{%
    \begin{tabular}[#1]{@{}c@{}}#2\end{tabular}}

\usepackage[ruled,vlined]{algorithm2e}
\usepackage{multirow}
\usepackage{makecell}
\usepackage{boldline}
\usepackage{nicefrac}
\usepackage{booktabs}

\begin{document}


\preprint{To appear in Mathematics and Mechanics of Solids (\url{https://doi.org/10.1177/10812865211055504})}

\title{\Large{Predicting Peak Stresses In Microstructured Materials Using Convolutional Encoder-Decoder Learning}}

\author{Ankit Shrivastava}
    \affiliation{Department of Civil and Environmental Engineering, Carnegie Mellon University}

\author{Jingxiao Liu}
    \affiliation{Department of Civil and Environmental Engineering, Stanford University}

\author{Kaushik Dayal}
    \affiliation{Department of Civil and Environmental Engineering, Carnegie Mellon University}
    \affiliation{Center for Nonlinear Analysis, Department of Mathematical Sciences, Carnegie Mellon University}
    \affiliation{Department of Materials Science and Engineering, Carnegie Mellon University}

\author{Hae Young Noh}
    \email{noh@stanford.edu}
    \affiliation{Department of Civil and Environmental Engineering, Stanford University}
    \affiliation{Department of Civil and Environmental Engineering, Carnegie Mellon University}

\date{\today}


\begin{abstract}
	This work presents a machine learning approach to predict peak-stress clusters in heterogeneous polycrystalline materials.
	Prior work on using machine learning in the context of mechanics has largely focused on predicting the effective response and overall structure of stress fields. 
	However, their ability to predict peak stresses -- which are of critical importance to failure -- is unexplored, because the peak-stress clusters occupy a small spatial volume relative to the entire domain, and hence requires computationally expensive training.
	This work develops a deep-learning based Convolutional Encoder-Decoder method that focuses on predicting peak-stress clusters, specifically on the size and other characteristics of the clusters in the framework of heterogeneous linear elasticity.
	This method is based on convolutional filters that model local spatial relations between microstructures and stress fields using spatially weighted averaging operations. 
	The model is first trained against linear elastic calculations of stress under applied macroscopic strain in synthetically-generated microstructures, which serves as the ground truth.
	The trained model is then applied to predict the stress field given a (synthetically-generated) microstructure and then to detect peak-stress clusters within the predicted stress field.
	The accuracy of the peak-stress predictions is analyzed using the cosine similarity metric and by comparing the geometric characteristics of the peak-stress clusters against the ground-truth calculations.
	It is observed that the model is able to learn and predict the geometric details of the peak-stress clusters and, in particular, performed better for higher (normalized) values of the peak stress as compared to lower values of the peak stress.
	These comparisons showed that the proposed method is well-suited to predict the characteristics of peak-stress clusters. 
\end{abstract}

\maketitle

\section{Introduction}\label{intro}

This paper aims to use a deep learning-based method to predict peak-stress clusters in heterogeneous linear elastic materials.
Specifically, we consider a simplified model of a polycrystalline microstructure as being composed of many grains, and each grain has an anisotropic linear elastic response.
The response of each grain is identical, except that the crystal lattice -- and hence the elastic stiffness tensor -- in each grain has a rotation with respect to some reference grain.
In this idealization, the polycrystalline microstructure can be described by a piecewise-constant field $\bfR(\bfx)$ where $\bfR\in SO(2)$ is a rotation that maps the lattice in the reference grain to the orientation of the local lattice at the material point located at $\bfx \in \mathcal{A}$, where $\mathcal{A} \subset \mathbb{R}^2$ is the domain.
Then, given the elastic stiffness tensor $\bfC_0$ of the reference, the elastic stiffness at a point $\bfx$ is simply the standard tensor transformation of $\bfC_0$ through $\bfR$.
This provides a piecewise-constant field $\bfC(\bfx)$ as the heterogeneous stiffness field. This idealization is standard in homogenization, e.g., \cite{liu2004homogenization}.

The heterogeneity in the material properties generically leads to stress fluctuations and is the subject of numerous works in homogenization, e.g. \cite{lipton2016uncertain}.
The stress fluctuations can lead to stress concentrations and regions of high stress, and it is critical to predict the peak stresses that are the dominant drivers of material failure.
While numerical methods for elasticity, such as the finite element method \cite{bathe2006finite} or Fourier-based schemes \citep{de2017finite, zeman2017finite, Lebensohn2012,peng2020effective} can solve for the stress field -- and the peak stresses as subsequent post-processing -- these methods are relatively expensive, particularly if applied to a large number of microstructures drawn from a random distribution.
Therefore, machine learning-based methods have recently been proposed as an alternative that can potentially be less expensive, particularly the process of prediction after training \citep{mianroodi2021teaching, mao2020designing, yang2018microstructural, lubbers2017inferring, li2018transfer,mozaffar2019deep}. 
In these methods, the machine learning models\footnote{In machine learning, the statistical expression that expresses a relationship between input variable and response variable is referred to as a ``model''. 
	It is often distinguished from the word ``method'', which describes the entire process of preprocessing the input,  making a prediction using the model given the input, and post-processing the model prediction for further application.} 
were trained on the data (``ground truth'') obtained from numerical methods to learn the statistical relationship between the microstructure and the corresponding response to applied loads.
When used for prediction, these trained models are much faster than the existing numerical schemes. 

However, these proposed methods have not been characterized in terms of their ability to predict the peak stresses, as opposed to the effective response that is simply the relation between the applied load and the overall deformation.
For instance, an interesting recent contribution is \cite{frankel2020prediction} that uses a Convolutional Long Short Term Memory (LSTM) based deep-learning model to predict stress fields in microstructures and then used this to predict the effective stress response (spatially averaged stress) under applied strain.
Although the model shows remarkable fidelity and potential, the size of the microstructure input is limited in the number of grains, making it difficult to explore the question of peak stress prediction in a statistical sense. 
Moreover, training a similar Convolutional LSTM model to predict the stress fields on large-sized microstructures will require an extremely large dataset, making it very expensive.
Further, an important next step is to characterize the ability of learning-based methods to not only predict averaged response correctly but the peak stresses as well, which are the dominant drivers of material failure.

While the ability of machine learning-based methods to predict peak stresses is largely unexplored, it has recently started to receive attention.
An important recent contribution is \cite{Mangal2018AppliedMaterials}, wherein they proposed a supervised learning-based random forest model to classify grains in a given microstructure as high stress or not.
The input is in the form of specified features such as the number of neighbors, Schmid factors, and so on. 
Although this approach has shown promising results, it relies on decisions by the user and can suffer from user bias. 
Hence, we propose an unsupervised learning-based approach to this problem.

Specifically, we propose an approach to predict peak-stress clusters in a given microstructure using a combination of a deep learning-based Convolutional Encoder-Decoder (CED) model and an image processing algorithm.
The design of the CED model used to predict stress fields is inspired by Convolutional Autoencoder.
The Autoencoder is an unsupervised deep learning algorithm that aims to reconstruct a given input. 
Convolutional Autoencoder (CAE) is a specific type of Autoencoder that uses convolutional filters to reconstruct their input. 
The input to the Convolutional Autoencoders is a single channel (such as grayscale) or multichannel image (such as RGB images)\footnote{In this work, we refer to a 3D array of size $H \times W\times C$ as a multichannel image, where $C$ is the number of channels. 
	The input to Convolutional Encoder Decoder model is a four channel image with size $H=128$, $W=128$ and $C=4$. 
	Similarly ground-truth for the model is a single-channel image with size $H=32$, $W=32$ and $C=1$.
	While discussing dimension $H \times W$, we refer to them as image resolution.}\label{fn:channel}. 
Since most of the images' features are locally-spatially correlated, the convolutional filters can extract the local features, which a deep fully-connected layer cannot capture. Due to the ability of CAE to learn image-image mapping with arbitrary accuracy, it is heavily used in image compression, de-noising, generation, and segmentation \cite{ronneberger2015u, badrinarayanan2017segnet}. 

Typical polycrystalline microstructures under an applied load show peak-stress clusters that are sparse and are governed by the local spatial features of the microstructure as suggested by \cite{Rollett2010StressPolycrystals}. 
For a 2D elasticity problem, a given microstructure and corresponding stress field (von Mises) can be treated as images. 
Hence architectures based on convolution operation are appropriate. 
Therefore, to find a mapping between the microstructure and corresponding stress field, CAE models are used. 
As the output and input are different for our problem, we refer to our model as Convolutional Encoder-Decoder.

First, the trained CED model is used to predict the stress field for a given microstructure.
Then a cluster detection algorithm is used to detect peak-stress clusters within the predicted stress field.  
The mean squared error was used to analyze the accuracy and to critically investigate scenarios where the model makes highly inaccurate predictions. 
Finally, the accuracy of the peak-stress clusters' geometric structure was analyzed using the cosine similarity metric and by separately comparing the clusters' size, shape, orientation, and location in the learning-based predictions against the ground-truth numerical elasticity calculations.
Visualization of the filters of the first convolutional layer were used to understand the essential elements of the microstructure as identified by the machine learning model.

\paragraph*{Organization.}

In Section \ref{method}, the basic principles used to design the CED model and algorithm to detect peak-stress clusters are discussed. 
Section \ref{experiment} describes the procedure to generate the required dataset and the principle behind designing the architecture of CED. 
Section \ref{results} discusses the evaluation of our model. 
In Section \ref{Feature_visualisation}, the role of filters in the first layer of the CED model is discussed. 
Concluding remarks are provided in Section \ref{conclusion}.

\section{Convolutional Encoder-Decoder Based  Method to Predict Peak-Stress Clusters}\label{method}

\begin{figure*}[htp]
	\centering
	\includegraphics[width=\textwidth]{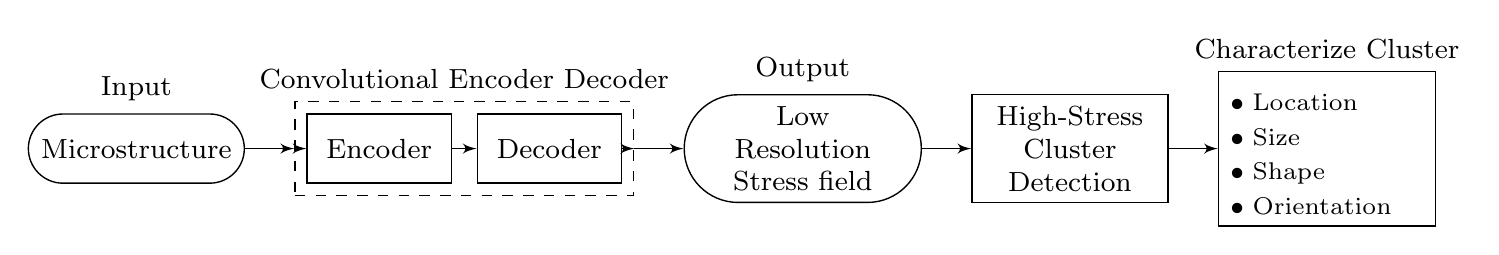}
	\caption{Flow chart of the proposed method to predict peak-stress clusters in a given microstructure.  
		In this work, we use the term ``model'' to refer to the Convolutional Encoder-Decoder and the word ``method'' to refer to the complete procedure from predicting stress fields from microstructures using CED to detecting and characterizing peak-stress clusters inside the predicted stress fields.}\label{flowchart}
\end{figure*}

The peak-stress clusters in a given microstructure are predicted using the approach shown in Figure \ref{flowchart}. 
First, the trained CED model predicts a low-resolution von Mises stress field from a given microstructure. 
Then model predictions are fed into a cluster detection algorithm to detect peak-stress clusters. 
The detected peak-stress clusters are then characterized by their location, shape, size, and orientation. 
The training and designing of the model are explained in later sections.

\subsection{Input and Ground-Truth for the Model}\label{sec:input-gt}

\begin{figure*}[htp]
	\centering
	\includegraphics[width=\textwidth]{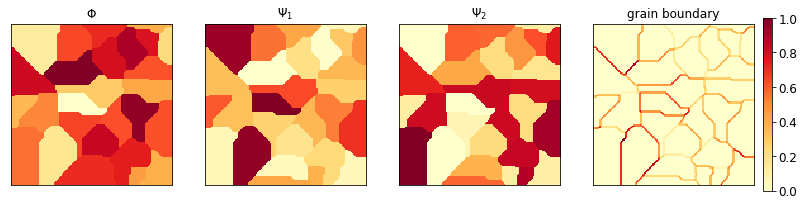}
	\caption{Images of resolution $128 \times 128$ representing Euler angles $(\Phi,\Psi_1,\Psi_2)$ and the grain boundary structure. 
		These images are scaled between $[0,1]$ before being fed into the model.}
	\label{input}
\end{figure*}

The CED model's input is a microstructure of 3D polycrystals in a 2D domain, i.e., we have $3$ Euler angles at every spatial location $\bfx$.
The input therefore contains information about the orientation and grain boundaries. 
The ground-truth is a low-resolution von Mises stress field. 
Since a large number of samples are required to train the deep learning model and test its validity, we create a synthetic dataset with known properties. 
The Dream3D package \cite{Groeber2014} was used to generate synthetic single-phase microstructures.
The package uses a set of predefined statistics for grain geometry and Euler angle distributions to generate these microstructures; see Table \ref{tab:stats} for more details. 
A microstructure generated  by the package is a three-channel image\footnotemark[\value{footnote}], each channel representing one of the Euler angles, $(\Phi,\Psi_1,\Psi_2)$. 
The grain boundary structure is a single channel image obtained using the crystal orientation information. 
Each point of the grain-boundary image is the $L_2$-norm of the gradient of Euler angles at that point.
The three Euler angle images and the grain-boundary image (Figure \ref{input}) are used as the input to the CED model. 
The primary purpose of including the grain-boundary image was to improve the CED model accuracy. 
The extra information about the grain boundaries as an input to the model became necessary as the model could not deduce this accurately by itself.

\begin{figure*}[htp]
	\centering
	\subfloat[ \label{fig:von_mises}]{\includegraphics[width=.35\textwidth]{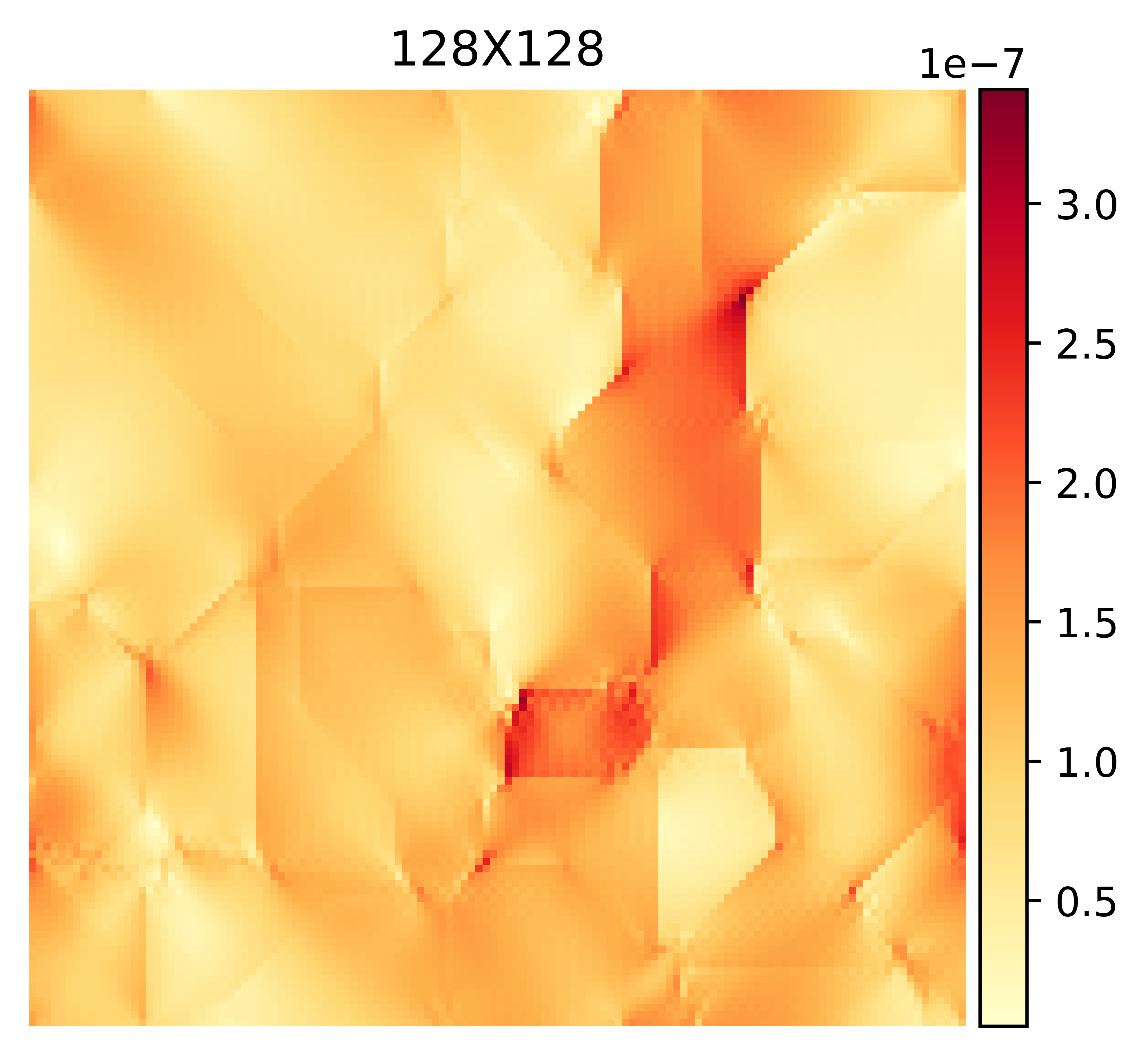}}
	\hspace{5em}
	\subfloat[ \label{fig:von_mises_low}]{\includegraphics[width=.35\textwidth]{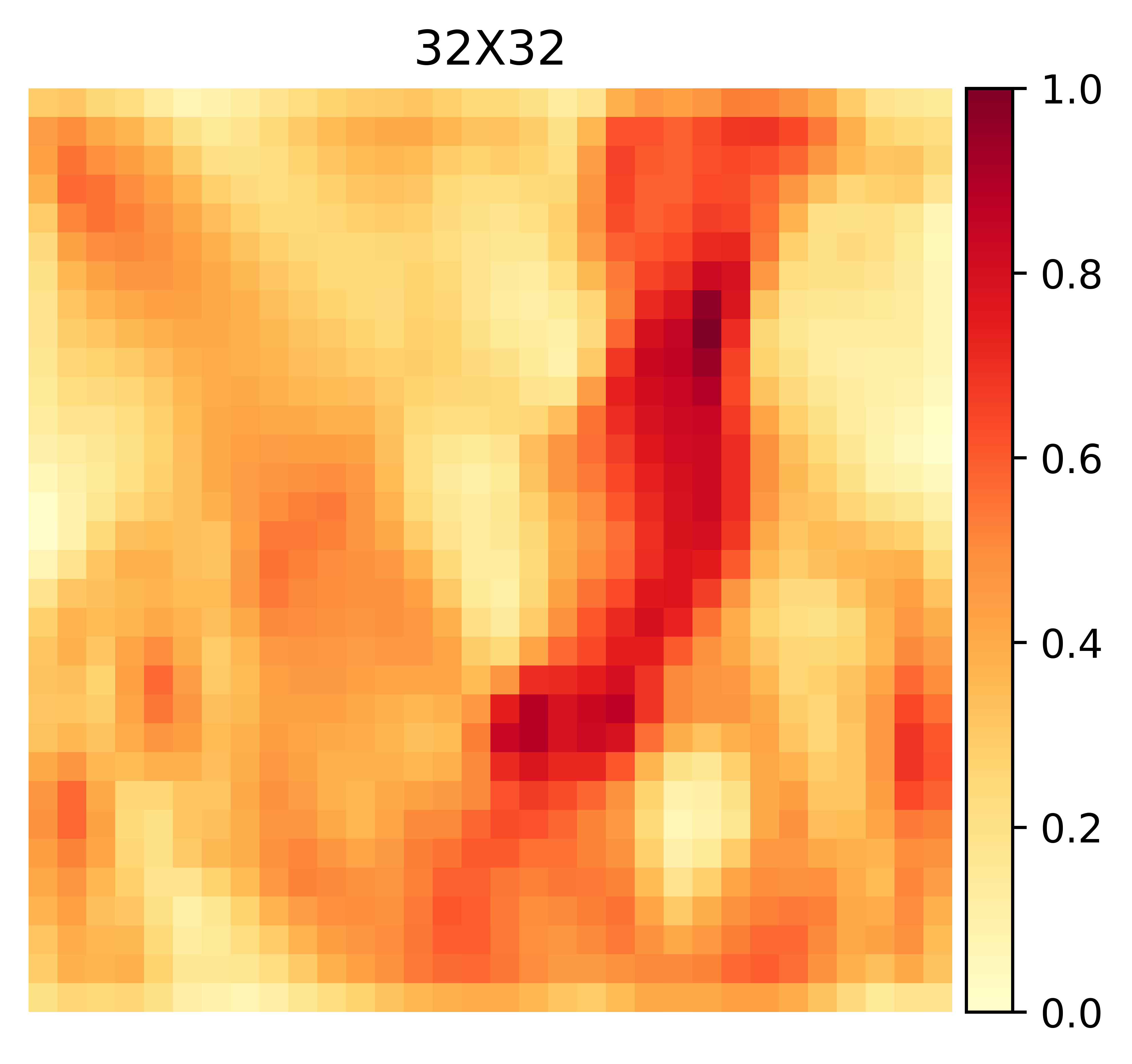}}
	\caption{(a) von Mises stress field and (b) low-resolution von Mises stress field obtained by piecewise-averaging used as ground truth. This stress field is scaled between $[0,1]$.}
	\label{output}
\end{figure*}

The stress field is computed using linear elasticity for each of the microstructures under an applied macroscopic strain, and the stress field is used to calculate the von Mises stress field. 
The von Mises stress fields are single-channel images of the same resolution as the microstructures. 
However, training any deep learning model with large-sized images will require a vast amount of training data that requires significant computational resources and memory. 
Hence, we use piecewise-averaged stress fields, which approximates the von Mises stress fields, as shown in Figure \ref{output}. 
These low-resolution stress fields are the ground-truth for the model. 
The piecewise averaging of stress fields preserves the relative character of the peak-stress clusters while making the model training efficient. 
The CED model architecture is designed considering the size of the microstructures and the ground-truth stress fields.

\subsection{Design Criteria for Convolutional Encoder Decoder}\label{sec:Model}

The CED model contains two main components, Encoder and Decoder. 
Both components consist of a sequence of layers with multiple filters that operate on its input and subsequently pass its output to the next layer. 
The output from each layer is referred to as the feature map. 
Details of the operations used in these components are explained in the later sections.

\begin{figure}
	\includegraphics[width=.95\textwidth]{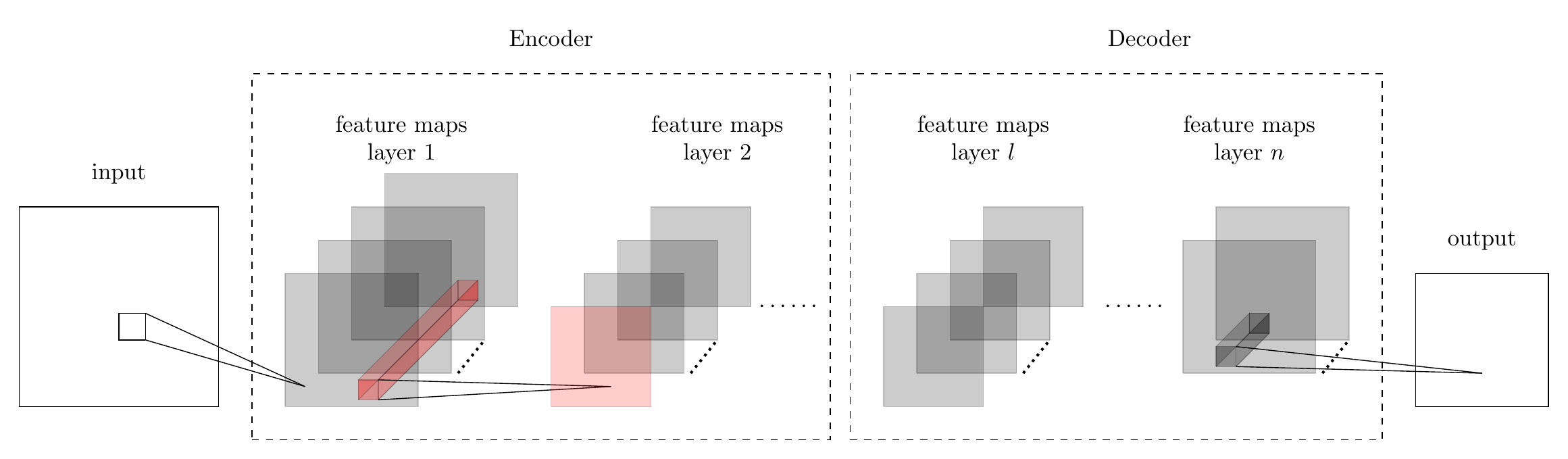}
	\caption{A general framework of the proposed CED model. 
		The input layer consists of a four-channel image representing a microstructure as shown in Figure \ref{input}. 
		The output layer consists of a single channel image representing a low-resolution von Mises stress field as shown in Figure \ref{fig:von_mises_low}. 
		The Encoder consists of multiple alternating convolutional and max-pooling layers. 
		The Decoder consists of multiple transpose convolution layers.}
	\label{CED}
\end{figure}

\subsubsection{Encoder}\label{sec:encoder}

The Encoder uses a series of convolutional and max-pooling layers to capture spatial information such as grain boundary, Euler angles, and grain junctions in the neighborhood of a point in the microstructure. 
It then converts the information into low-resolution with high-description feature maps. 
A convolutional layer performs convolution operation using multiple filters followed by a nonlinear activation. 
These filters separately extract the high-level features from the layer input by convolving (weighted dot product) it with a sliding window, as shown in Equation \ref{eq:convolution}. 
The hyperparameters, such as the number of filters in the layer, size of each filter, padding $P$ on input, and stride $S$, are predefined. 
Their choice is explained in Section \ref{hyperparameter}.

Since the convolutional operations are linear, their combination will make the model capture only the linear relationship between the microstructures and the low-resolution stress fields, ignoring any nonlinear relations. 
Hence, outputs from each convolution operation are fed into a nonlinear activation function. 
These activation functions make the layer nonlinear, further assisting in extracting the nonlinear features from its input \cite{Nair2010}. 
The output size from the nonlinear activation function is the same as the input size to the function.

\begin{align}
	n \times n \times K_{in} &\equiv \text{kernel size of each convolutional filter} \nonumber\\
	n_{in} \times n_{in} \times K_{in} &\equiv \text{Dimension of input to the layer, } \mathbf{C\_{in}} \nonumber\\
	n_{out} &= \frac{n_{in}-n+2P}{S}+1 \label{eq:conv_out}\\ 
	K_{out} &\equiv \text{number of filters} \nonumber\\
	n_{out} \times n_{out} \times K_{out} &\equiv \text{Dimension of output after convolution, } \mathbf{z} \nonumber\\
	n_{out} \times n_{out} \times K_{out} &\equiv \text{Dimension of output after operating $\relu$, } \mathbf{C\_{out}} \nonumber\\
	\mathbf{z}_{ijk} &=   \sum _{c = 1} ^{K_{in}}\sum _{b = 1} ^{n}\sum _{a = 1} ^{n} \mathbf{w}_{abck} \mathbf{C\_{in}}_{(i-a)(j-b)c} + \mathbf{w}_{k} \quad \forall \ i,j \in [1,n_{out}] \quad \text{and} \quad k \in [1,K_{out}] \label{eq:convolution}\\
	\mathbf{w}_{abck} \text{ and } \mathbf{w}_k & \text{ are the weights and bias respectively of $k^{th}$ filter in the layer.} \nonumber
\end{align}
\begin{align}\label{eq:relu}
	\mathbf{C\_{out}}_{ijk} = \relu(\mathbf{z}_{ijk}) = \max(0,\mathbf{z}_{ijk})  \quad \forall \ i,j \in [1,n_{out}] \quad \text{and} \quad k \in [1,K_{out}]
\end{align}

For the CED model, $\relu$ activation (Equation \ref{eq:relu}) is used as non linear activation. 
Unlike other activation functions such as $\sigmoid$ and Hyperbolic tangent function $\tanh$, $\relu$ functions do not suffer from the vanishing gradient problem \cite{Goodfellow-et-al-2016}. 
$\relu$ is also much more computationally efficient \cite{qiu2018frelu} making the back-propagation computationally efficient. 
This combination of convolution and nonlinear activation operations summarizes the information of input ($\mathbf{C\_{in}}$) by detecting its features in the form of feature maps ($\mathbf{C\_{out}}$).

In the Encoder, the number of filters in the subsequent layers is increased as more and more layers are added. 
Since the size of the feature maps is decreasing, further reducing the amount of the microstructure information, the number of filters is increased to compensate for the reduction in feature size. 
However, using a large number of convolutional filters in deeper layers will increase the number of parameters in the model, which can cause over-fitting. 
Maxpool layer is often used instead of a convolutional layer to avoid a large number of parameters. 
The maxpool operation, given by Equation \ref{eq:maxpool}, separately operates on each of its input $\mathbf{M\_{in}}$ using a sliding window, extracting the most important information while keeping the number of model parameters constant.
This operation further reduces the size of its input while keeping important information. 
During training, maxpool also prevents fluctuations in values of parameters $ \mathbf{w}$ in convolutional layers due to a slight change in its input. 
This further helps in achieving stable convergence while training.

\begin{align}
	n_{out} \times n_{out} \times K_{out} &\equiv \text{Dimension of its input, } \mathbf{M\_{in}} \nonumber\\
	m \times m \times 1 &\equiv \text{kernel size of the maxpool filter} \nonumber\\
	m_{out} &= \frac{n_{out}-m+2P}{S}+1 \label{eq:max_out}\\
	\mathbf{M\_{out}}_{ijk} &= \max\{ \mathbf{M\_{in}}_{(i-a)(j-b)k} \mid a \in [1,m] \text{,} b \in [1,m] \}  \label{eq:maxpool}\\ &\quad \forall \ i,j \in [1,m_{out}] \quad \text{and} \quad k \in [1,C_{out}] \nonumber
\end{align}

\subsubsection{Decoder}\label{sec:decoder}

The Decoder converts feature maps obtained from the Encoder into low-resolution von Mises stress fields using a series of transpose convolution layers, each performing transpose convolution operation given by Equation \ref{eq:tconvolution}. 
The transpose convolution operation decompresses the information from its input using weighted interpolation. 
For the same reason as explained in Section \ref{sec:encoder}, the transpose convolution operation is followed by a nonlinear activation. 
$\relu$ is used as the activation function for all layers except for the last layer.

\begin{align}
	n \times n \times C_{in} &\equiv \text{kernel size of each transpose convolutional filter} \nonumber\\
	n_{in} \times n_{in} \times C_{in} &\equiv \text{Dimension of layer's input, } \mathbf{A} \nonumber\\
	n_{out} &= (n_{in}-1) \times S + n - 2P \label{eq:tconv_out}\\
	\mathbf{z}^{t}_{ijk} &=   \sum _{c = 1} ^{C_{in}}\sum _{b = 1} ^{n}\sum _{a = 1} ^{n} \mathbf{w}^{t}_{abck} \mathbf{A}_{(i+a)(j+b)c} + \mathbf{w}^{t}_{k} \quad \forall \ i,j \in [1,n_{out}] \quad \text{and} \quad k \in [1,C_{out}] \label{eq:tconvolution}\\
	\mathbf{w}^{t}_{abck} \text{ and } \mathbf{w}^{t}_k & \text{ are the weights and bias respectively of $k^{th}$ filter in the layer.} \nonumber
\end{align}

In the last transpose convolution layer, the transpose convolution operation is followed by a $\tanh$ activation function given by Equation \ref{eq:tanh}. 
The model does not suffer from a vanishing gradient by including $\tanh$. 
The $\tanh$ function is used in the last layer of the model since it is more nonlinear than the $\relu$. 
This helps in closely resembling the ground-truth with just a few layers, whereas the $\relu$ being piecewise linear will require more layers to resemble the same ground-truth. 
The $\tanh$ is preferred over the $\sigmoid$ as the former has a stronger gradient than the latter, which reduces the vanishing gradient problem while training. 
The output from the last layer is the prediction of a low-resolution stress field.

\begin{align}\label{eq:tanh}
	\tanh(\mathbf{z}_{ijk}) = \frac{e^{\mathbf{z}_{ijk}} - e^{-\mathbf{z}_{ijk}}}{e^{\mathbf{z}_{ijk}} + e^{-\mathbf{z}_{ijk}}}
\end{align} 

\subsubsection{Training Model Parameters}\label{sec:OB}

The Mean Squared Error ($\MSE$) function is considered as the objective function for training and validating the model.  
The $\MSE$ represents how close the model predictions are to the ground-truth for a given input.
First, the complete dataset is divided into training, validation, and test sets.
The $\MSE$ computed over the training set, the validation set, and the test set are referred to as training error, validation error, and test error, respectively. 
The goal while training the model is to find optimum parameters, $\mathbf{w^*}$ of the CED that minimizes the training error, see Equation \ref{eq:min}. 
The training input data $\mathbf{X}$ represents a set of microstructures, and the ground-truth $\mathbf{Y}$ represents a set of low-dimensional stress fields for the corresponding input.

\begin{align}
	\mathbf{w}^* &= \arg\min_{\mathbf{w}}\MSE(\mathbf{w}) \label{eq:min} \\
	\MSE(\mathbf{w}) &= \underset{i \in [1,N]}{\mathop{\mathbb{E}}}\left[\|\mathbf{CED}(\textbf{X}_{i},\mathbf{w}) - \mathbf{Y}_{n}\|^2\right] \label{eq:mse} \\
	\mathbf{CED}(\mathbf{X}_{i},\mathbf{w}) &\equiv \text{ the model prediction for a given input $\mathbf{X}_{i}$} \nonumber\\
	\mathbf{Y}_n &\equiv \text{corresponding ground-truth} \nonumber\\
	N &\equiv \text{ number of samples} \nonumber
\end{align}

The $\MSE$ is minimized using the mini-batch gradient descent scheme {\cite{Goodfellow-et-al-2016}}. 
Under this scheme, the training dataset is divided into sub-samples, called a batch, with a predefined batch size. 
The $\MSE$ over a batch is used to compute the gradients and update the model parameters. 
Compared to gradient descent, the mini-batch gradient descent helps avoid local minima due to the noisy update, further helping achieve faster convergence. 
It also avoids the problem of GPU memory restrictions as suffered by standard gradient descent. 
The mini-batch gradient is chosen over the stochastic gradient descent to achieve smoother training and avoid expensive latency in CPU-GPU data transfer during every update.
As compared to the stochastic gradient descent, the mini-batch gradient descent has a higher likelihood to fall into a local minimum {\cite{Goodfellow-et-al-2016}}. 
However, this problem can be controlled by tuning the learning rate and the batch size. These parameters help control the amount of noise added during the gradient update, which helps avoid local minima.
For the gradient computations, we have used the ADAM scheme {\cite{kingma2014adam}}. 
Compared to other gradient computation schemes, ADAM provides faster convergence by adapting the geometry of the objective function. 

During the training, the validation and training errors are evaluated regularly to analyze if the model is underfitting or overfitting.
Model hyperparameters are tuned based on underfitting and overfitting as explained in Section \ref{hyperparameter}
Overfitting is also avoided using an early-stopping mechanism whereby the model training is stopped as soon as validation error increases.

\subsubsection{Hyperparameter Selection for the Model}\label{hyperparameter}

We use the behavior of the validation and training errors during the training to select the model hyperparameters. The different values of the model hyperparameters will lead to a different model architecture. 

The hyperparameters such as the number of convolutional, maxpool, transpose convolutional layers and the number of filters in each layer are governed by the validation error. 
A large number of layers in the Encoder can lead to a significant loss in the information about the input as each layer compresses the information. 
A large number of layers in the Decoder can lead to the inclusion of noise during the weighted interpolation. 
In both cases, the validation error increases.  
A small number of layers in both the Encoder and the Decoder can cause both the training error and the validation error to increase (underfitting).

Similarly, in each layer, a small number of filters will have a small number of parameters, which will increase the chance of the model getting stuck in a local minimum causing high validation error. 
A large number of filters can either lead to an increase in validation error (overfitting) or no significant reduction in validation error while making training expensive.

The kernel sizes of the filters in all the layers are governed by the validation error, the input dimension, and the ground-truth dimension.
The stride and the padding of convolutional filters, maxpool filters, and transpose convolutional filters are governed by Equation {\ref{eq:conv_out}},  Equation {\ref{eq:maxpool}}, and Equation {\ref{eq:tconv_out}}, respectively.
The stride is always less than or equal to kernel size in convolutional layers to ensure that information from every pixel is considered while making predictions. 

Hence, optimum values of the hyperparameters are chosen to achieve a minimum validation error.
Finally, the architecture that gives the minimum validation error during the training is used to predict the stress fields.

\subsection{Detection of Peak-Stress Clusters}\label{sec:hotspots_detect}

\begin{figure}
	\includegraphics[width=\textwidth]{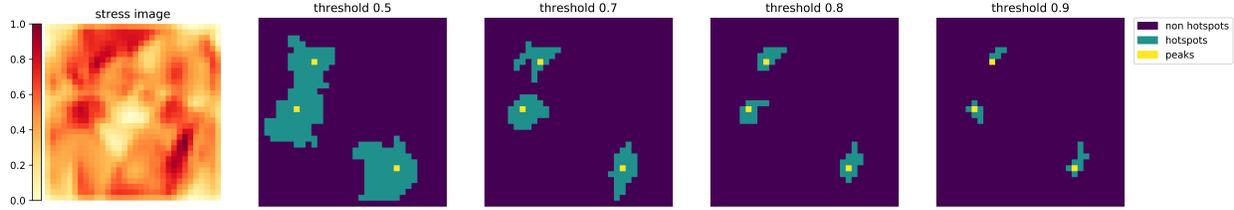}
	\caption{Stress clusters around three different peaks in a predicted stress field for different thresholds.}
	\label{fig:hotspot}
\end{figure}

The peak-stress clusters are detected inside the stress fields predicted by the trained CED model. 
We characterize the peak-stress clusters as regions around a peak in the stress fields with values above a certain threshold. 
The low-resolution stress fields obtained by averaging von Mises stress fields, see Section {\ref{sec:input-gt}}, and the predicted stress fields from the model are not smooth, and there are small fluctuations.
The stress fields are smoothed using a Gaussian filter to avoid small fluctuations. 
If the smoothing is not performed, small fluctuations in the stress fields will be detected as peaks, which would lead to the erroneous assignment of these small fluctuations as the actual peak.
Further, the algorithm used to detect peak-stress clusters around stress peaks will be computationally expensive.

After the smoothing, the stress fields are fed into a peak detection algorithm which then returns the location of multiple peaks inside stress fields.
For every peak, the clusters around it are obtained by assigning the value of $1$ if the pixel value is higher than the predefined threshold times the corresponding peak stress. 
The result is a binary image with multiple disconnected regions representing stress clusters. 
The stress clusters are then separately labeled using the connected component labeling algorithm as shown in Figure \ref{fig:hotspot}. 

\subsection{Characterizing Stress Clusters}\label{sec:hotspots_characterize}

\begin{figure}
	\centering
	\includegraphics[width=.75\textwidth]{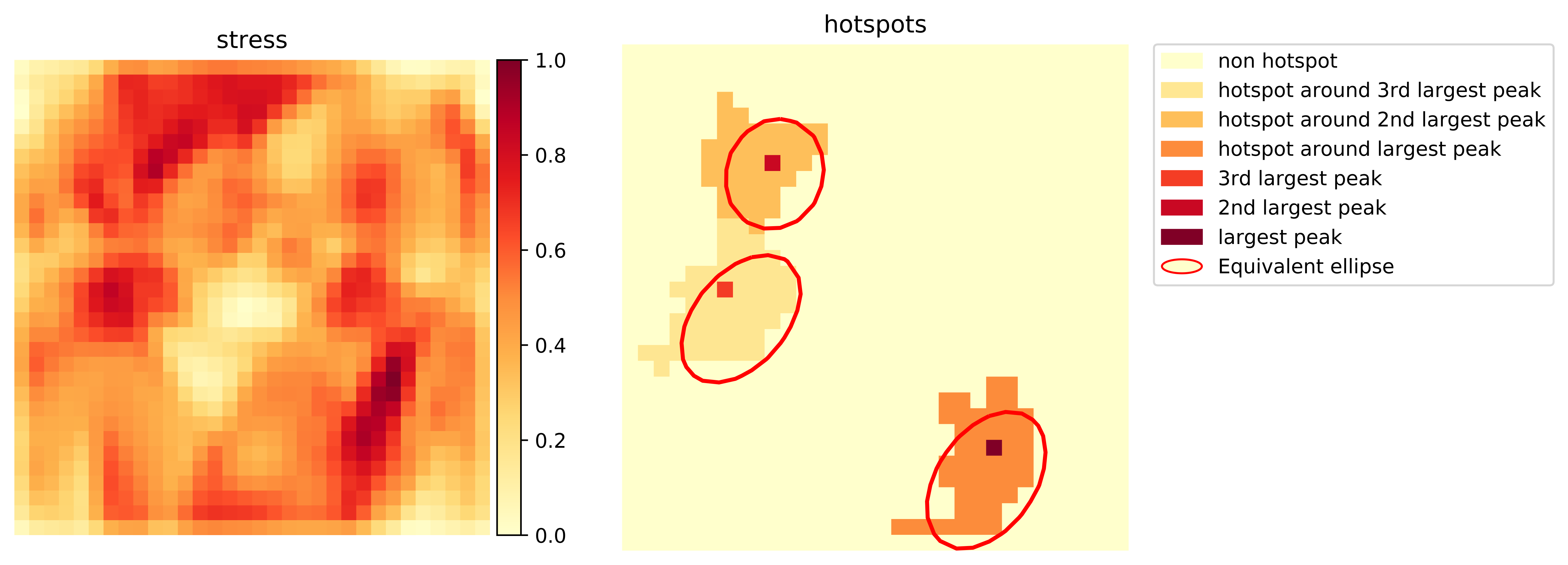}
	\caption{Equivalent ellipse of the stress clusters corresponding to each peak.}
	\label{fig:ellipse}
\end{figure}

The peak-stress clusters corresponding to each peak are characterized by their location, size, shape, and orientation. 
Since failure initiates at the peak-stress clusters, detailed information of such clusters is important. 
This characterization provides geometric details of peak-stress clusters inside a stress field. 
The location of the corresponding peak characterizes the location of a stress cluster. 
The ratio of the number of pixels in the cluster to the total number of pixels of the stress field characterizes the cluster size. 
The equivalent ellipse of the peak-stress cluster is obtained by the least-square fitting; see {\cite{rosin1993note}} for more details. 
Equation {\ref{eq:lsq}} is minimized to obtain the coefficients $\mathbf{z}_i$ in Equation {\ref{eq:ellipse}}.
The ellipse is centered around the mean of the $n$ pixels with coordinates $X$ and $Y$ inside the cluster. 
The ellipse's orientation, $\theta$ and axis ratio, AR are then obtained using the coefficients $\mathbf{z}$ as shown in Equation {\ref{eq:theta}} and {\ref{eq:ar}}. 
It is important to note that the right-hand side of Equations {\ref{eq:ellipse}} is set to $1$ in order to avoid a trivial solution.
Hence, the value of $\mathbf{z}$ that is obtained represents the original coefficients of the ellipse scaled by a constant.

\begin{align}
	\mathbf{z}_1 X^2 + \mathbf{z}_2 XY &+ \mathbf{z}_3 Y^2 + \mathbf{z}_4 X + \mathbf{z}_5 Y = 1 \label{eq:ellipse} \\
	\arg\min_{\mathbf{z}} &\|\mathbf{A}\mathbf{z} - \mathbf{b}\|_2^2 \label{eq:lsq} \\ 
	\text{where, } \mathbf{z}&=[\mathbf{z}_1, \mathbf{z}_2, \mathbf{z}_3, \mathbf{z}_4, \mathbf{z}_5]^T \nonumber \\
	\mathbf{b}&=[1,1,1,1,1]^T \nonumber \\
	\mathbf{A_i} &= [X_i^2 , X_i Y_i, Y_i^2, X_i, Y_i] \quad \forall i \in [1,n]\nonumber \\
	n &\equiv \text{ Number of points inside a cluster}& \nonumber \\
	\theta &=  \tan^{-1}\left(\frac{\mathbf{z}_3-\mathbf{z}_1-\sqrt{(\mathbf{z}_1-\mathbf{z}_3)^2 + \mathbf{z}_2^2}}{\mathbf{z}_2}\right) \label{eq:theta}\\
	\text{AR} &= \frac{\mathbf{z}_1+\mathbf{z}_3+\sqrt{(\mathbf{z}_1-\mathbf{z}_3)^2 + \mathbf{z}_2^2}}{\mathbf{z}_1+\mathbf{z}_3-\sqrt{(\mathbf{z}_1-\mathbf{z}_3)^2 + \mathbf{z}_2^2}} \label{eq:ar}
\end{align}

\section{Numerical Experiments}\label{experiment}

\subsection{Generating the Data Set}

A dataset of $10000$ samples containing the input and the ground-truth are generated as explained in Sections \ref{sec:gen_in} and \ref{sec:gen_gt}. 
The dataset is further divided into $7000$ for training, $1000$ for validation, and $2000$ for testing. 
Before being used for training or prediction, the dataset is pre-processed.

Initially, the model was trained using a dataset with $40000$ samples. 
However, the test error using $40000$ samples was $0.0188618$, whereas, for $10000$ samples, it was  $0.0197408$. 
Hence, increasing the data did not improve the accuracy significantly but increased the training time by a factor of four, which acted as a bottleneck to the tuning of model architectures. 
It also required a large amount of memory to store the data and results for post-prediction analysis. 
Hence, in the later sections, the model training, the hyperparameter tuning, and the evaluation are performed using the dataset with $10000$ samples.

\subsubsection{Generating Input for the Model}\label{sec:gen_in}

For the current problem, microstructures of size $128 \times 128 \times 3$ were generated using Dream3D with parameters in Table \ref{tab:stats}. 
Many microstructures generated with the dimensions of $64$ and $32$ were similar in terms of the grain distribution. 
The dimension of $128$ ensured that the generated microstructures are sufficiently different from each other. 
Next, the information about the grain boundaries inside a microstructure was obtained using the $L_2$-norm of the gradient of the Euler angles for each pixel inside the microstructure. 
The grain boundary information, along with the microstructure, was used as an input. 
The final dimension of the input was $128 \times 128 \times 4$ (Figure \ref{input}).

\begin{table}[htp]
	\centering
	\begin{tabular}{|cc|cc|}
		\Xhline{2\arrayrulewidth}
		\textbf{Preset Statistic} : & Primary Equiaxed  & \textbf{Mu} : & 2.3   \\ \hline
		\textbf{Max cutoff} : & 4 & \textbf{Sigma} : & 0.4   \\ \hline
		\textbf{Bin step size} : & 1 & \textbf{Min/Max cutoff}  : & 4 \\ 
		\Xhline{2\arrayrulewidth}
	\end{tabular}
	\caption{Values in the table represent Dream3D StatsGenerator parameters to generate microstructures. 
		The Primary-Equiaxed parameter represents a single-phase material with grains of shape ellipsoid having almost the same length in all directions. The equivalent sphere diameter of grains in microstructures is randomly selected from a truncated lognormal distribution defined using Mu, Sigma, and Min/Max cutoff. 
		The Bin step size further segregates the lognormal distributions into multiple bins, where each bin will have a corresponding shape distribution with default values set by the Preset Statistic parameter. 
		Finally, the Euler angles in all three directions are randomly generated from a uniform distribution.}
	\label{tab:stats}
\end{table}

\subsubsection{Generating ground-truth for the Model}\label{sec:gen_gt}

\begin{figure*}[htp]
	\centering
	\includegraphics[width=.3\textwidth]{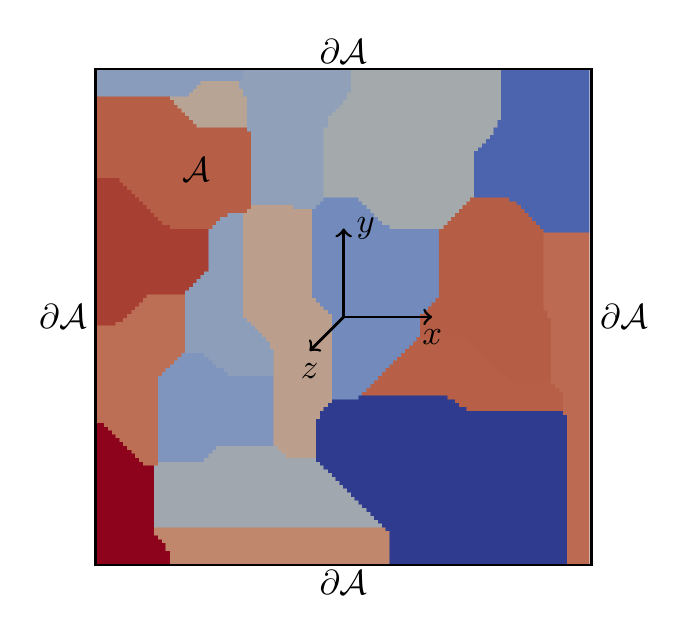}
	\caption{A typical microstructure sample in the XY plane. $\mathcal{A}$ and $\partial \mathcal{A}$ represent domain and boundary of the microstructure.}
	\label{fig:applied_load}
\end{figure*}

Linear elasticity (Equation {\ref{eq:le}}) was applied to generate the von Mises stress fields using a Fourier method {\cite{Lebensohn2012}}.
The non-dimensionalized fourth-order elasticity tensor considered represents copper with cubic symmetry (Table {\ref{tab:elastic}}).

\begin{equation}
	\begin{aligned}
		\nabla \cdot \boldsymbol{\sigma}(\mathbf{x}) &= 0  &&\forall \mathbf{x} \in \mathcal{A} \\
		\boldsymbol{\sigma}(\mathbf{x}) &= \mathbf{C}(\mathbf{x}) : \boldsymbol{\epsilon}(\mathbf{x}) &&\forall \mathbf{x} \in \mathcal{A} \\
		\boldsymbol{\epsilon}(\mathbf{x}) &= \mathbf{e}(\mathbf{u^*(\mathbf{x})}) + \mathbf{E} &&\forall \mathbf{x} \in \mathcal{A} \\
		\mathbf{e}(\mathbf{u^*(\mathbf{x})}) &= \frac{\nabla \mathbf{u}^*(\mathbf{x}) + \nabla \mathbf{u}^*(\mathbf{x})^T}{2} \quad &&\forall \mathbf{x} \in \mathcal{A} \\ 
		\mathbf{u}^*(\mathbf{x})& &&\text{Periodic in } \mathcal{A}\\
		\mathbf{E} &=  
		\begin{bmatrix}%
			0 & 0 & 0 \\
			0 & 0 & 0 \\
			0 & 0 & \num{1e-4} \\
		\end{bmatrix}
		\label{eq:le} 
	\end{aligned}
\end{equation}
where $\mathbf{E}$ is the applied macroscopic strain, $\mathbf{u^*(\mathbf{x})}$ is the periodic fluctuation displacement field, and $\boldsymbol{\epsilon(\mathbf{x}})$ is the strain field.

\begin{table}[htp]
	\centering
	\begin{tabular}{|c|c|c|}
		\Xhline{2\arrayrulewidth}
		$c_{1111}$& $c_{1122}$ & $c_{2323}$  \\ \hline
		$1.0000$ & $0.7209$ & $0.4477$ \\ 
		\Xhline{2\arrayrulewidth}
	\end{tabular}
	\caption{Normalized values of the elastic stiffness tensor components for cubic symmetry.} 
	\label{tab:elastic}
\end{table}

The size in pixels of the von Mises stress field is $128\times128\times1$ (Figure \ref{fig:von_mises}). 
For the reasons explained in Section \ref{sec:input-gt}, a moving average filter of size $4 \times 4$ with stride $4$ and no padding was operated on the von Mises stress field to obtain a stress field of size $32 \times 32 \times 1$ (Figure \ref{fig:von_mises_low}).
This piecewise averaged stress field is the ground-truth for the CED model. 
Models that reconstructed the ground-truth of size $64 \times 64 \times 1$ and $128 \times 128 \times 1$ had a much higher reconstruction error as compared to  ground-truth of size $32 \times 32 \times 1$. 

\subsubsection{Preprocessing Input and Ground-Truth}\label{sec:gen_preprocess}

Before using the microstructures as the input and the low-resolution von Mises stress fields as the ground-truth for training and prediction respectively, their values are scaled to lie between $[0,1]$. 
The scaling is performed for each image and each channel separately. 
Since the model is nonlinear, the effect of the linear scaling of the input and ground truth during the gradient update of the parameters will be nonlinear. 
Hence, the scaled input and ground truth ensure that the pixels with very large values do not completely dominate during the computation of weights. 
For example, in the unscaled low-resolution stress fields, the pixel values range $[1\times10^{-4},13.68]$. 
Hence the pixels with values close to $13.68$ will dominate the training compared to the pixel with values $1\times10^{-4}$. However, when scaled between $[0,1]$, the pixels with high values will not dominate as much as they would have during the unscaled data.

\subsection{Convolutional Encoder Decoder Architecture}\label{sec:MA}

\begin{figure}[htp]
	\centering
	\includegraphics[width=.60\textwidth]{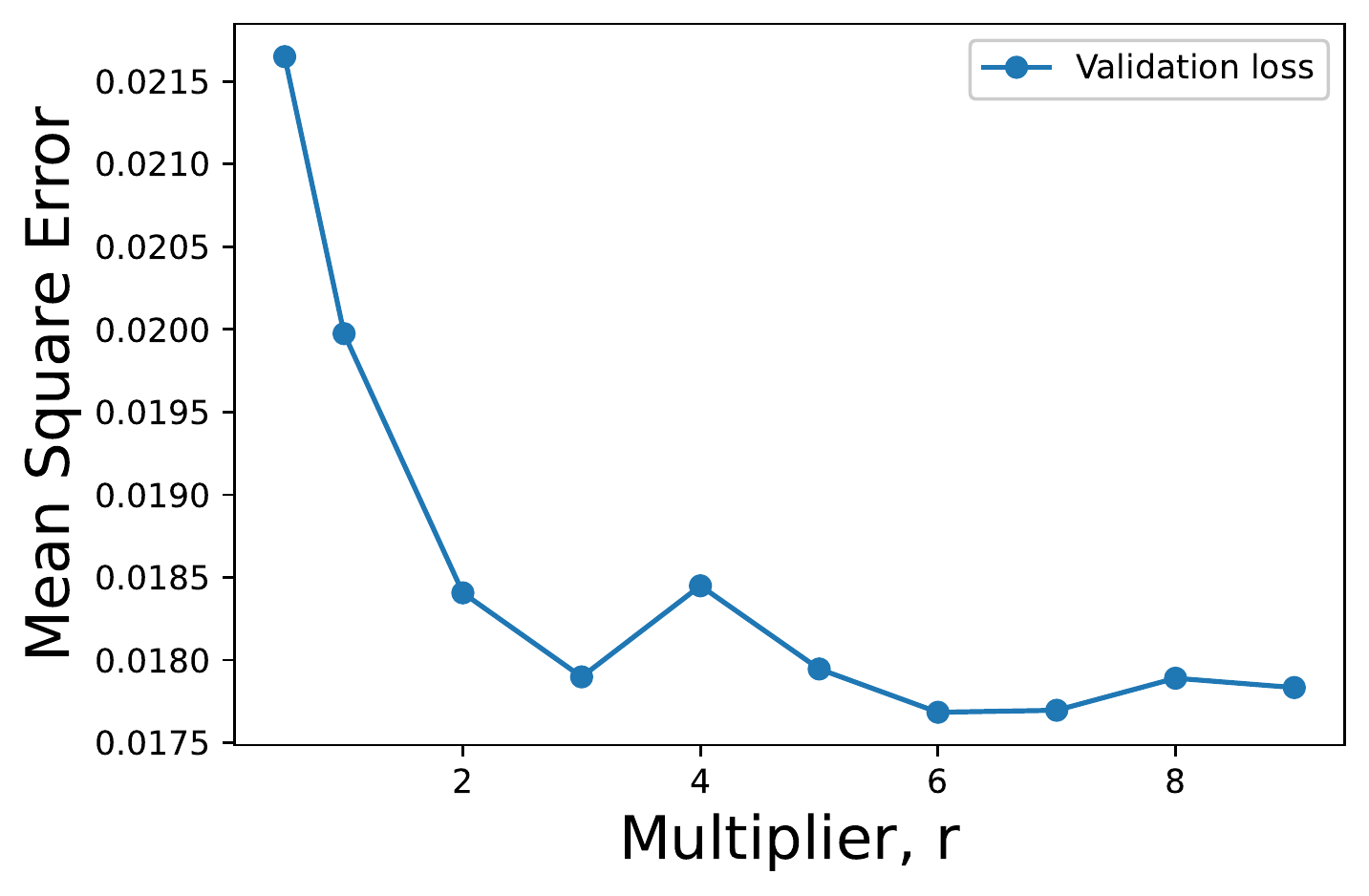}
	\caption{Validation loss for different values of the multiplier $r$. The number of filters is varied for each layer by multiplying $r$ with the number of filters in Table {\ref{tab:arch}}. The validation loss is for $r=3$ is $0.0178987$ and $r=9$ is $0.017834$. This suggests that the validation loss remains almost constant for $r>3$, hence $r=3$ is used as the optimal value. }
	\label{fig:channel_error}
\end{figure}

Various architectures of the CED model were designed based on the governing criteria discussed in Section \ref{sec:Model}.
The architectures were trained based on criteria discussed in Section \ref{sec:OB} and validation error was evaluated to tune the model hyperparameters.  
The hyperparameters used for training the architectures are shown in Table \ref{tab:optimisation}. 
Parameters such as the number of layers and kernel size of the filters were tuned based on the criteria discussed in Section \ref{hyperparameter}. 
The input of all layers was not padded as it implies fictitious information from outside the domain, causing the CED model to be less accurate. 
Without maxpool layers, the convergence of training loss was unstable; hence it was necessary to include it in the CED architecture. 
Once the number of layers, the function of each layer, and the kernel size of the filters were finalized, the number of filters for each layer, except for the output layer, was tuned by varying them as a multiple of variable $r$, see Table {\ref{tab:arch}}. 
The validation error against the variable $r$ is shown in Figure {\ref{fig:channel_error}}. 
With increase in $r$, the validation loss decreases till $r=3$ and then remains almost constant, i.e around 0.0178; hence $r=3$ is chosen as the optimal value. 
Since there is no significant difference in the performance between the model architecture for $r=1$ and $r=3$, the architecture with $r=1$ was chosen as the final CED architecture as it had a smaller number of filters in the first layer, making the analysis of the output of the first layer easier, see Section {\ref{Feature_visualisation}}.
PyTorch \cite{paszke2019pytorch} was used for designing and training the model. 

\begin{table}[H]
	\centering
	{
		\setlength\extrarowheight{2pt}
		\begin{tabular}{|c|c|c|c|c|c|}
			\hline
			\textbf{Components} &	\textbf{Operations in layers}	     & \textbf{Kernel size}      & \textbf{Stride} & \textbf{Number of filters for $r=1$} \\ 
			\hline
			\multirow{4}{*}{Encoder} 
			& \multicolumn{1}{c|}{Convolution + $\relu$}  & \multicolumn{1}{c|}{$8\times8\times4$}  & \multicolumn{1}{c|}{$8$} &  \multicolumn{1}{c|}{$64$} \\ 
			\cline{2-5}
			& \multicolumn{1}{c|}{Maxpool} 	 & \multicolumn{1}{c|}{$2\times2\times1$}  & \multicolumn{1}{c|}{$2$}    & \multicolumn{1}{c|}{$64$} \\ 
			\cline{2-5}
			& \multicolumn{1}{c|}{Convolution + $\relu$}  & \multicolumn{1}{c|}{$3\times3\times64$} & \multicolumn{1}{c|}{$1$}    & \multicolumn{1}{c|}{$16$} \\ 
			\cline{2-5}
			& \multicolumn{1}{c|}{Maxpool2D} 	 & \multicolumn{1}{c|}{$2\times2\times1$}  & \multicolumn{1}{c|}{$1$}    & \multicolumn{1}{c|}{$16$} \\ 
			\cline{1-5}
			\multirow{3}{*}{Decoder} 
			& \multicolumn{1}{c|}{\specialcell{Transpose Convolution  + $\relu$}}  & \multicolumn{1}{c|}{$2\times2\times16$}  & \multicolumn{1}{c|}{$2$}    & \multicolumn{1}{c|}{$64$} \\ 
			\cline{2-5}
			& \multicolumn{1}{c|}{\specialcell{Transpose Convolution  + $\relu$}} 	 & \multicolumn{1}{c|}{$3\times3\times64$}  & \multicolumn{1}{c|}{$1$}    & \multicolumn{1}{c|}{$8$} \\ 
			\cline{2-5}
			& \multicolumn{1}{c|}{\specialcell{Transpose Convolution  + $\tanh$}}  & \multicolumn{1}{c|}{$10\times10\times8$} & \multicolumn{1}{c|}{$2$}   & \multicolumn{1}{c|}{$1$} \\
			\cline{1-5}
		\end{tabular}
	}
	\caption{The architecture of the CED model. The input to all layers (including the model input) was not padded with any values as it reduced the model accuracy due to misleading information. The number of filters was tuned by varying it as a multiple of $r$. However, for evaluation and visualization purposes, $r=1$ is considered.}
	\label{tab:arch}
\end{table}

\begin{table}[H]
	\centering
	\begin{tabular}{|cc|cc|}
		\Xhline{2\arrayrulewidth}
		\hline
		\textbf{Batch size} : & 100  & \textbf{base learning rate} : & 0.0001   \\ \hline
		\textbf{Number of epochs} : & 4000 & \textbf{maximum learning rate} : & 0.1   \\ \hline
		\Xhline{2\arrayrulewidth}
	\end{tabular}
	\caption{Hyperparameters for optimization of $\MSE$. 
		One epoch means that each sample in the training set has been used once to update the model parameters. 
	}
	\label{tab:optimisation}
\end{table}

\section{Evaluation}\label{results}

The trained CED model with architecture as shown in Table \ref{tab:arch} was used to analyze 2000 separate test samples previously unseen during training. 
The performance of the peak-stress cluster prediction approach was analyzed using the MSE, cosine similarity metric, and by comparing clusters characteristics. 
First, we determined the microstructures from the test samples for which the model prediction is highly inaccurate. 
To do this, we characterized the microstructures in terms of their average grain diameter, average axis ratio, and average Euler angles. 
We then compared them against the MSE. 
Next, we excluded the microstructures with very high MSE while analyzing the accuracy inside the predicted stress fields. 
The cluster detection in such microstructures was not useful due to high inaccuracy in stress fields prediction. 
The accuracy of the clusters in stress fields was analyzed in two ways: first, by evaluating cosine similarity between rescaled predicted and ground-truth stress fields; and second, by comparing the characteristics of the peak-stress clusters detected in both prediction and ground-truth. 

\subsection{Error Analysis of the Predicted Stress Fields Using Mean Squared Error}

\begin{figure}[htp]
	\centering
	\includegraphics[width=.55\textwidth]{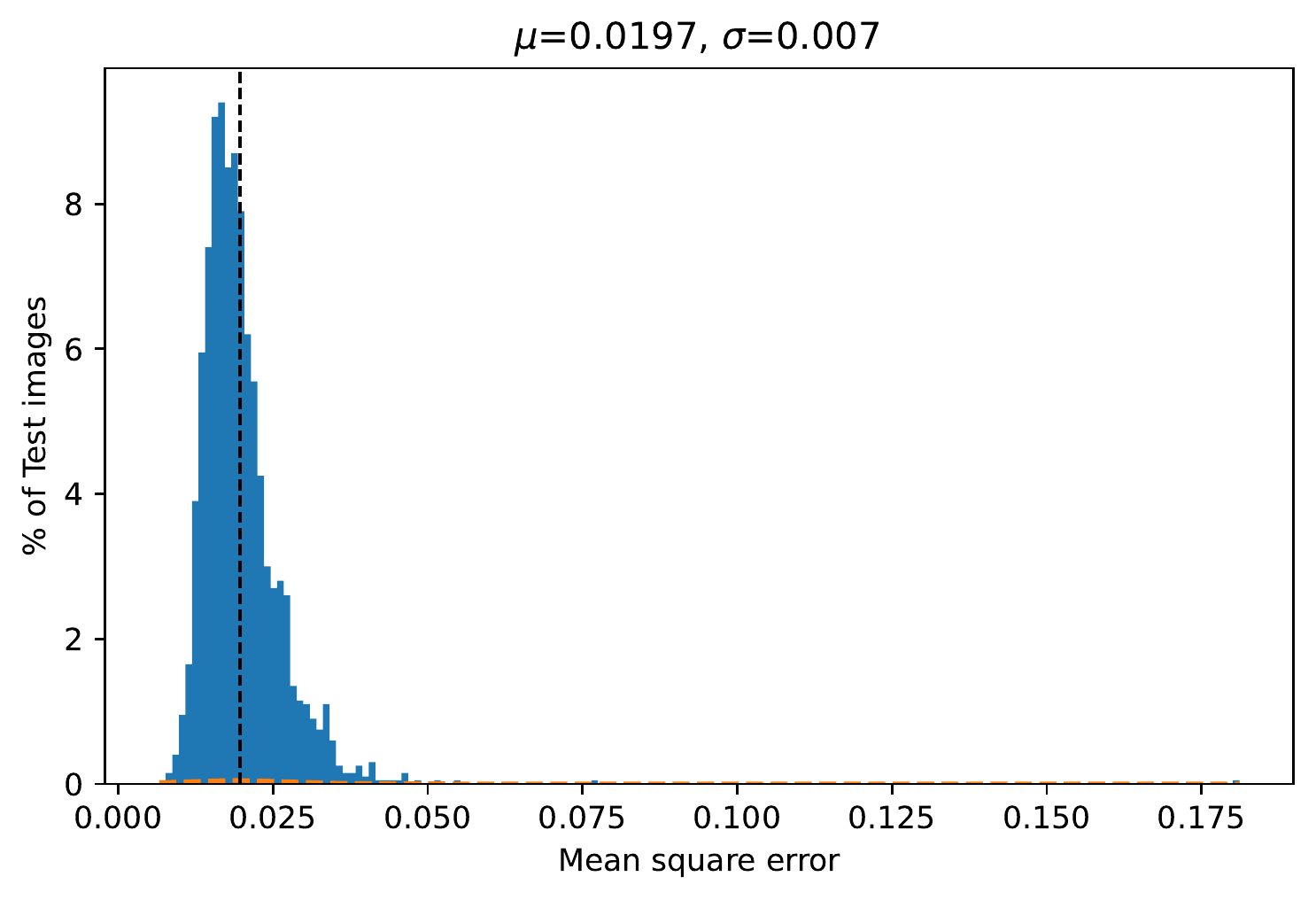}
	\caption{Histogram of MSE for $2000$ microstructure test samples. 
		The horizontal axis represents the MSE of each bin, and the vertical axis gives the percentage of samples for corresponding bins. 
		The parameters $\mu$ and $\sigma$ represent the mean and the standard deviation of the distribution, respectively. 
		Around $97.1\%$ of test samples have $\MSE$ below $0.0337$. 
		Hence, the average per-pixel error is around $\num{3.29e-05}$.
		The grain size inside the microstructures plays a major role in the accuracy of the model prediction.}
	\label{fig:mse}
\end{figure}

\begin{figure*}[htp]
	\centering
	\subfloat[Sample with MSE of $0.0066$. \label{fig:best_case}]{\includegraphics[width=.75\textwidth]{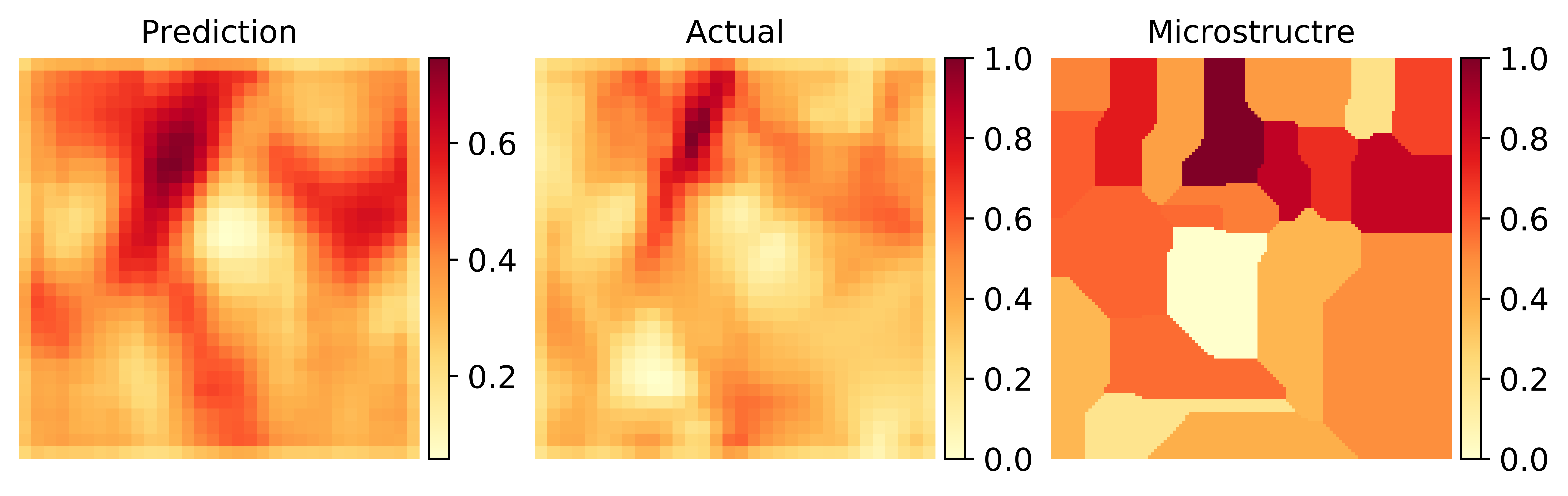}}
	\hfill
	\subfloat[Sample with MSE of $0.0204$.\label{fig:median_case}]{\includegraphics[width=.75\textwidth]{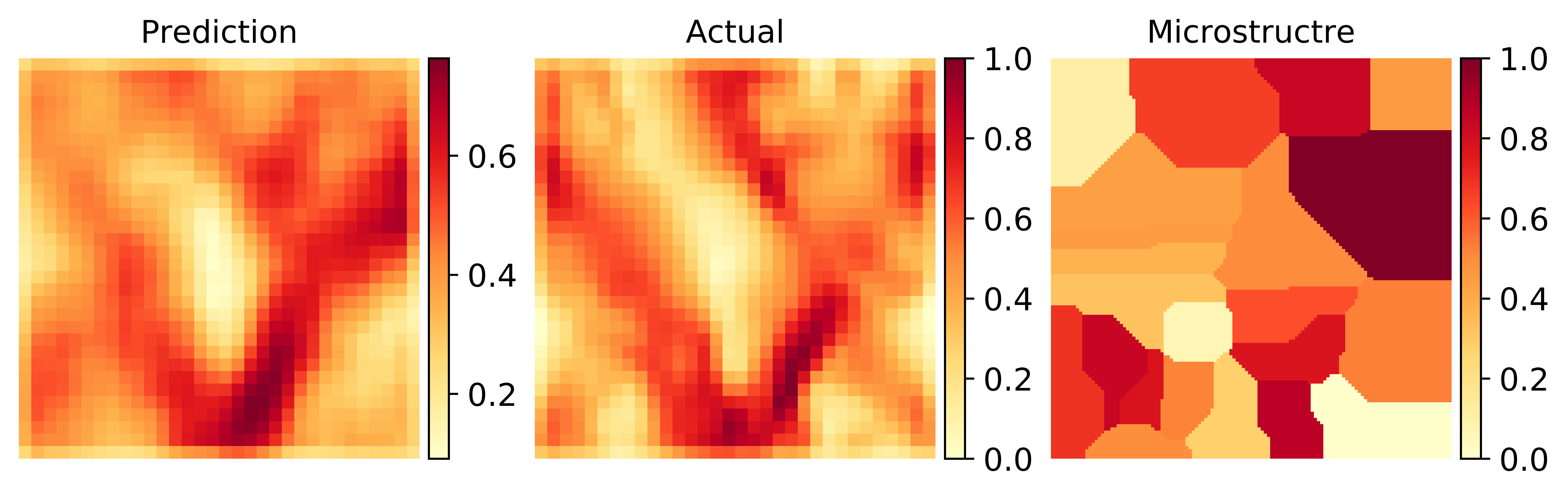}}
	\hfill
	\subfloat[Sample with MSE of $0.1811$\label{fig:worst_case}]{\includegraphics[width=.75\textwidth]{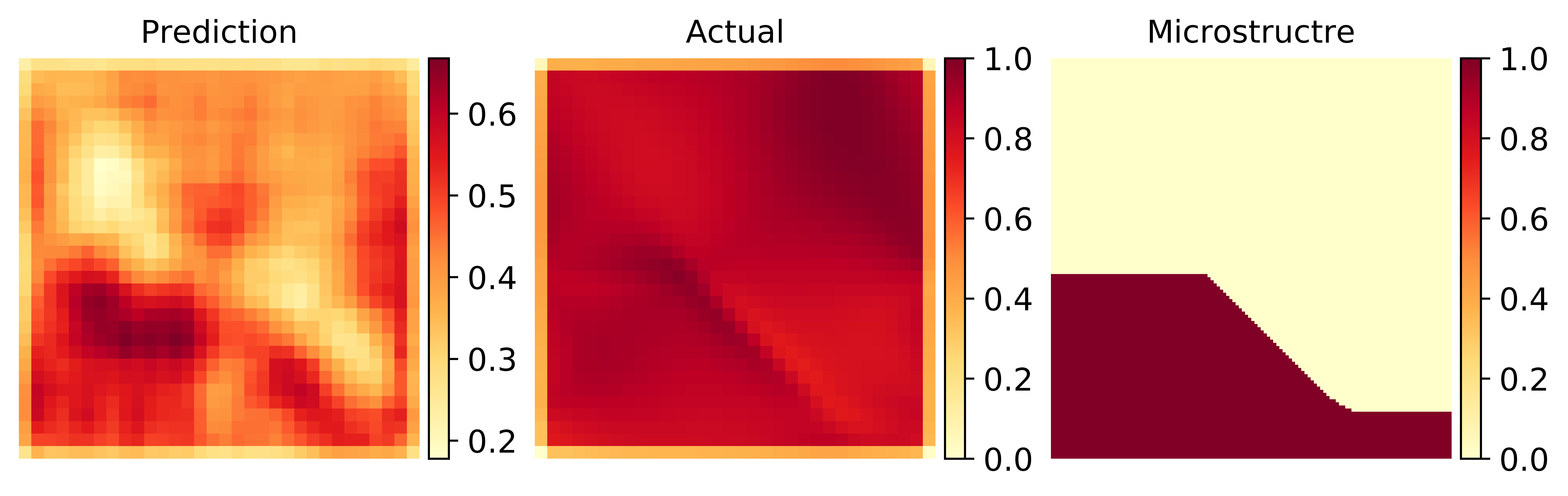}}
	\caption{The best (a), average (b), and worst (c) cases of low-resolution stress field prediction for the corresponding microstructure by the CED model. 
		The prediction suffers from high inaccuracy for microstructures with large grains (relative to the computational domain).}
	\label{fig:mse-visual}
\end{figure*}

\begin{figure*}[htp]
	\centering
	\subfloat[\label{fig:mse_grain_dia}]{\includegraphics[width=.45\textwidth]{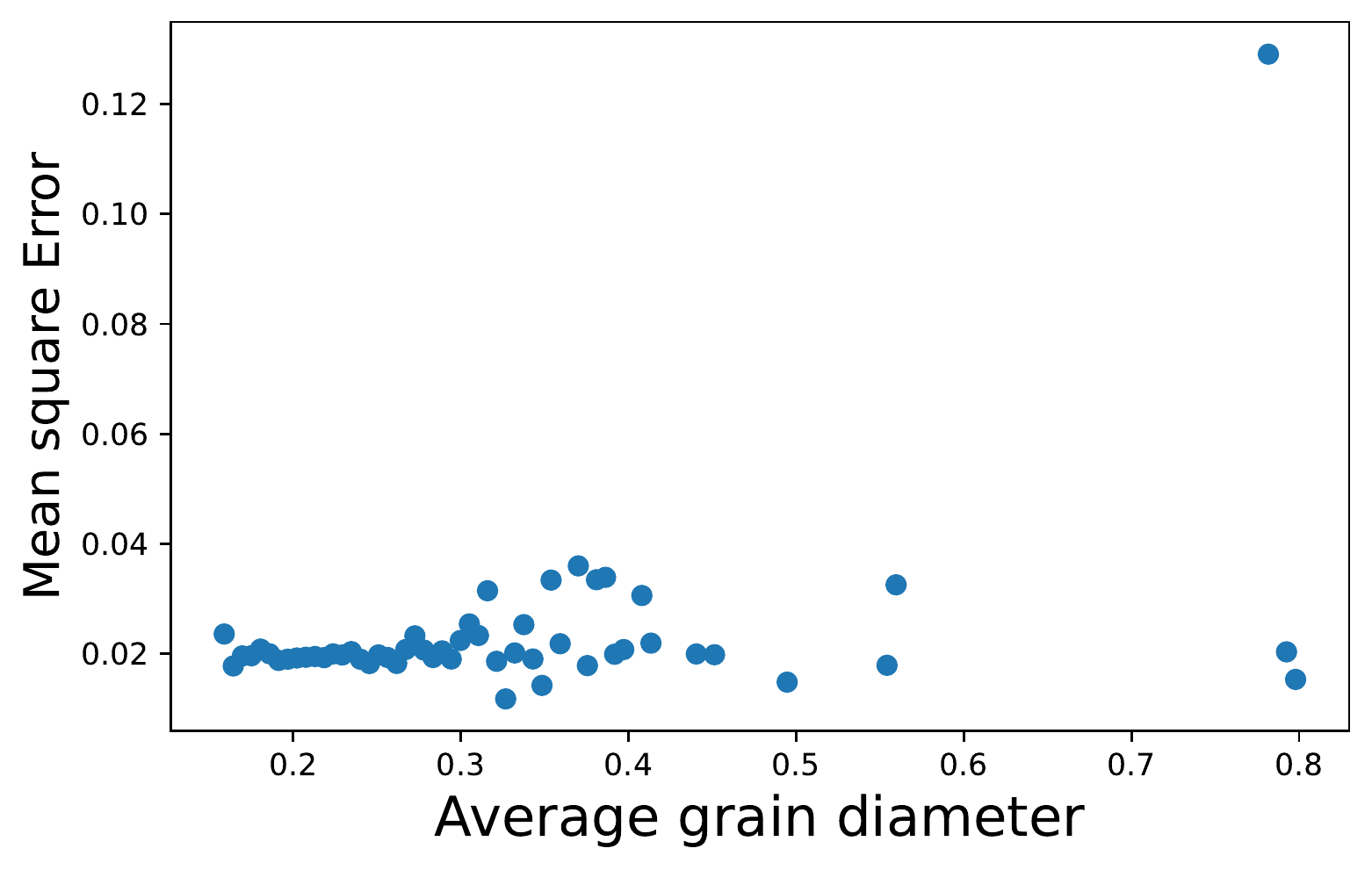}}
	\hfill
	\subfloat[\label{fig:mse_aspect}]{\includegraphics[width=.45\textwidth]{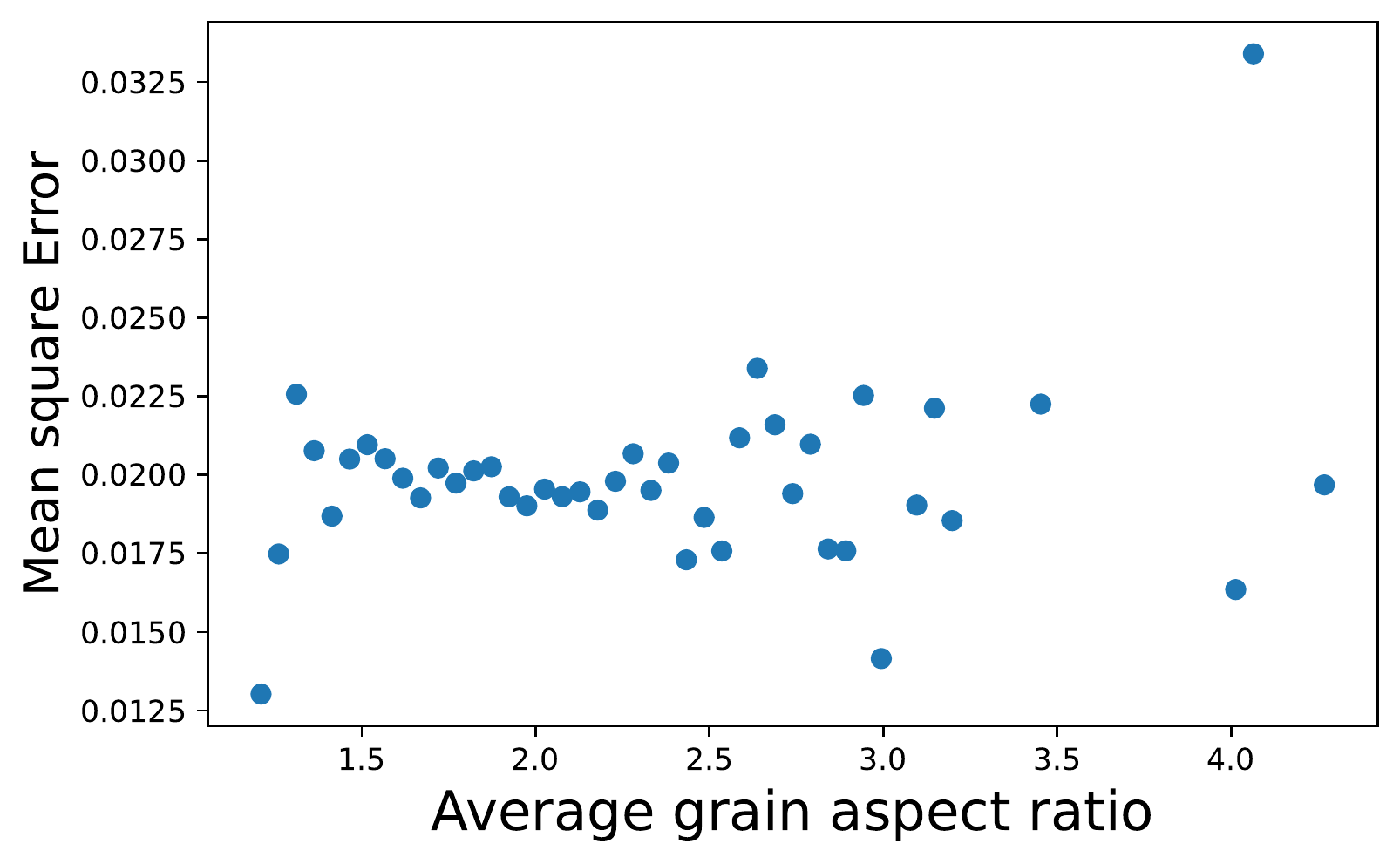}}
	\hfill
	\caption{Plot of average grain size and average aspect ratio against MSE. 
		The histogram of average grain size of the microstructures are computed. 
		For each bin of the histogram, the MSE and the average grain size of the microstructures belonging to the bin is collected and averaged. 
		The horizontal and vertical axes in (a) represents respectively the MSE and the average grain size (averaged over the histogram bins). 
		A similar calculation is performed to obtain (b).
		The grain diameter is normalized by the microstructure size, i.e., by $128$. 
		The MSE increases with mean grain size and mean grain aspect ratio, implying that the prediction was less accurate for microstructure with larger grains or elongated grains.}
	\label{fig:mse_ms}
\end{figure*}

\begin{figure*}[htp]
	\centering
	\subfloat[{Considering grains with diameter in range $[0, 1.12]$.\label{fig:mse_grain_dia_1}}]{\includegraphics[width=.45\textwidth]{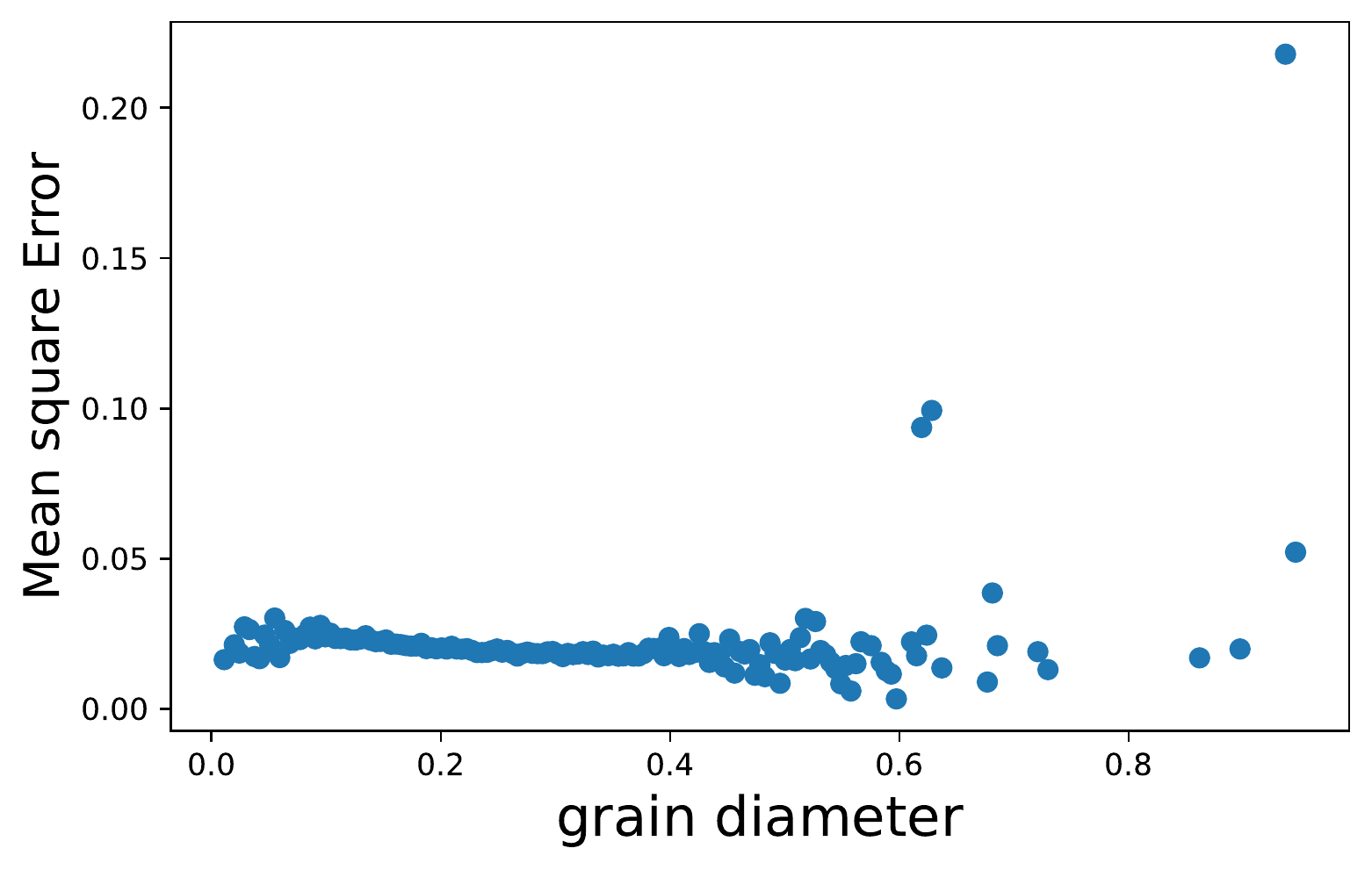}}
	\hfill
	\subfloat[{Considering grains with diameter in range $[0, 0.4]$.\label{fig:mse_grain_dia_2}}]{\includegraphics[width=.45\textwidth]{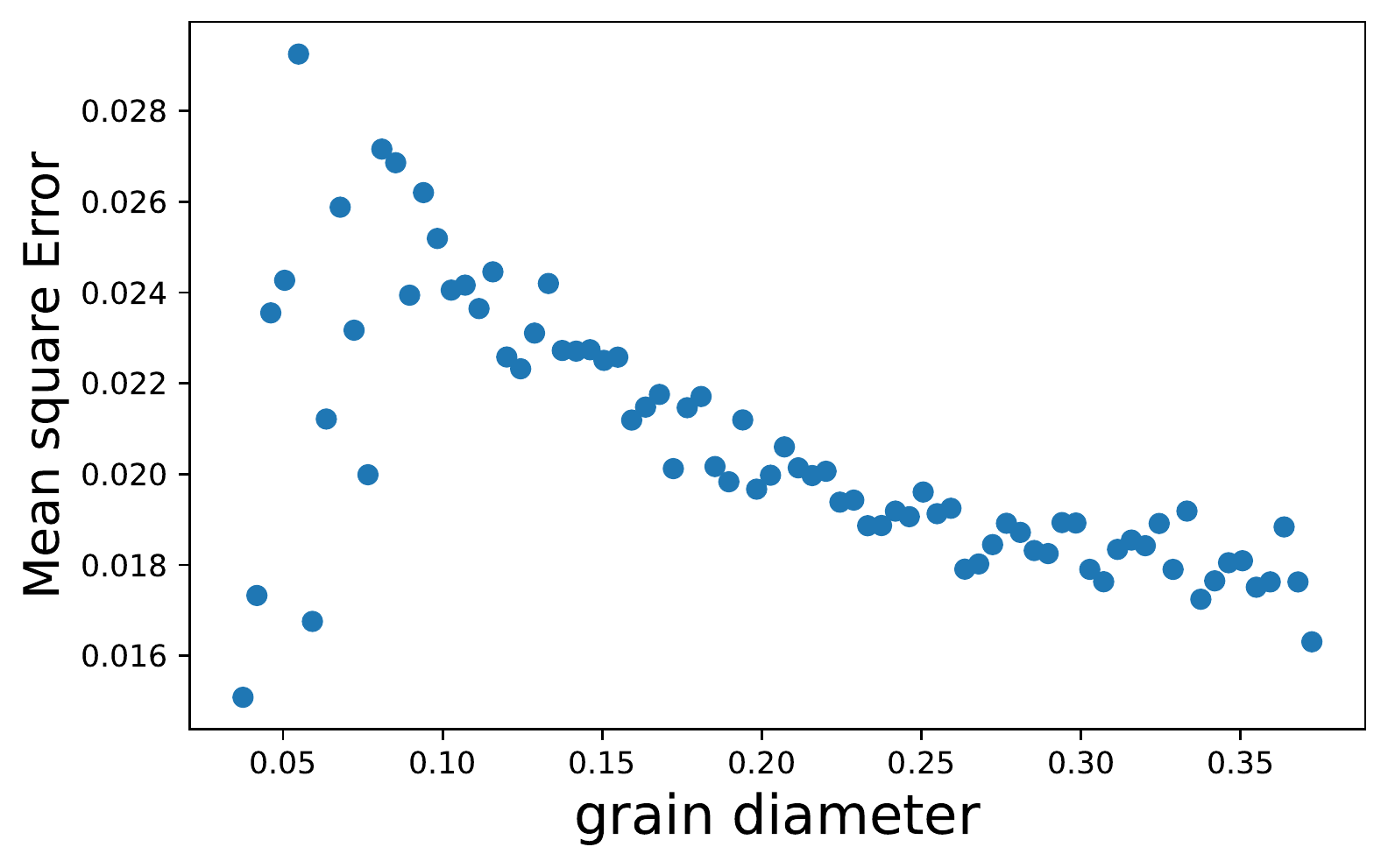}}
	\caption{Plot of grain diameter against MSE. 
		The plot is obtained in the same way as explained in Figure \ref{fig:mse_ms}.
		In (a), we consider all grains in the test set, and (b) considers only the subset of grains with diameter in the range $[0, 0.4]$. 
		(b) provides an understanding of the variation of MSE for smaller grains. 
		The MSE is high for very small grains and very large grains.}
	\label{fig:mse_grain}
\end{figure*}

The histogram of MSE for different microstructures is shown in Figure \ref{fig:mse}. 
The microstructures with average, largest, and smallest MSE and their corresponding model stress field prediction and ground-truth are shown in Figure  \ref{fig:mse-visual}. 
The MSE was studied against the microstructure's average grain diameter and average grain aspect ratio, as shown in Figure {\ref{fig:mse_ms}}.
For Figure {\ref{fig:mse_grain_dia}}, first the histogram of the average grain diameter is estimated.
The bins of the histogram are represented as a set $ \text{S} = \{\text{Bin}_1, \text{Bin}_2,
\text{Bin}_3, ...,
\text{Bin}_n\}$. 
Each bin in the histogram represents microstructures with similar average grain diameters. 
The average of MSE for all the microstructures that corresponds to $\text{Bin}_i$, $\forall \ i \in [1,n]$, is computed and represented as $ \text{MSE}_i$,  $\forall \ i \in [1,n]$.
These averaged MSE for each bin are plotted against the corresponding bin average grain diameter as shown in Figure {\ref{fig:mse_grain_dia}}.
The main reason to study the plots for each bin was to understand the effect of grains size on the model performance. If the grain size significantly impacts the model performance, it should be visible in the plots.
Figure  {\ref{fig:mse_grain_dia}} shows that the MSE increases with an increase in average grain diameter of the microstructure.
It suggests the prediction is better for microstructure with smaller grains. 
Similar behavior was observed when the MSE in each individual grain of the microstructures was analyzed against the corresponding grain diameter.
Figure \ref{fig:mse_grain_dia_1} shows that the error in the prediction of the stress distribution in larger grains was higher than in smaller grains. However, if the grain size is too small, then the error increases, as shown in Figure \ref{fig:mse_grain_dia_2}. 

The model performs better around medium-sized grains (size ranging from $[0.2,0.4]$) as grain junctions and grain boundaries provide sufficient feature information for prediction. 
In areas around larger grains, due to large homogeneous regions, the number of features available locally to the model is not enough to make predictions; hence, the model predicts stress with higher error. 
In the case of smaller grains, there are multiple sharp jumps in stress values around grains due to grain boundaries, which the model cannot capture effectively, leading to a higher error compared to medium-sized grains. 
This complete analysis suggests that if a microstructure consists mainly of medium-sized grains, then the model prediction error will be small.  
This behavior can also be observed by visual inspection, as shown in Figure \ref{fig:mse-visual}. 

\begin{figure*}[htp]
	\centering
	\subfloat[]{\includegraphics[width=.3\textwidth]{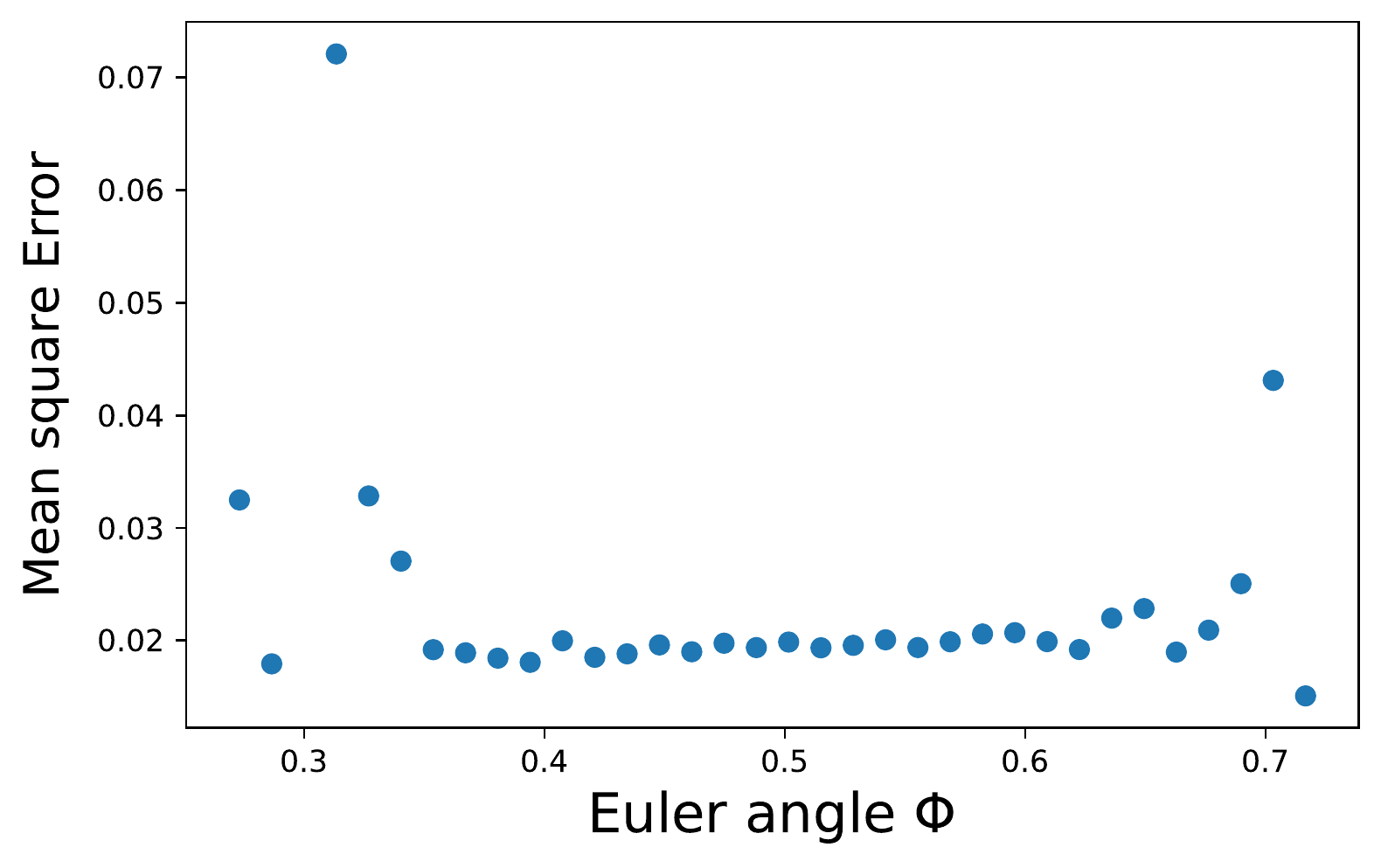}}
	\hfill
	\subfloat[]{\includegraphics[width=.3\textwidth]{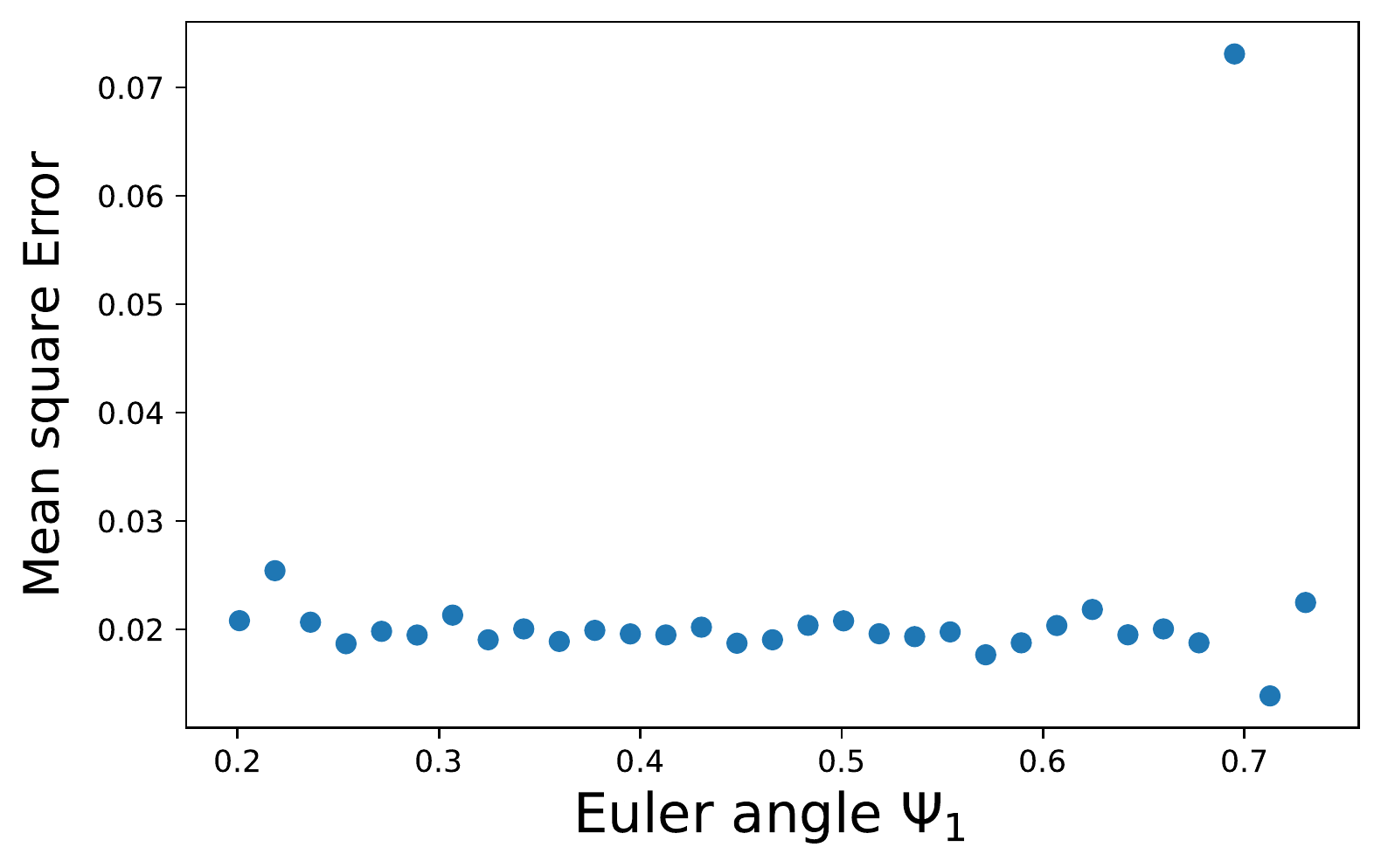}}
	\hfill
	\subfloat[]{\includegraphics[width=.3\textwidth]{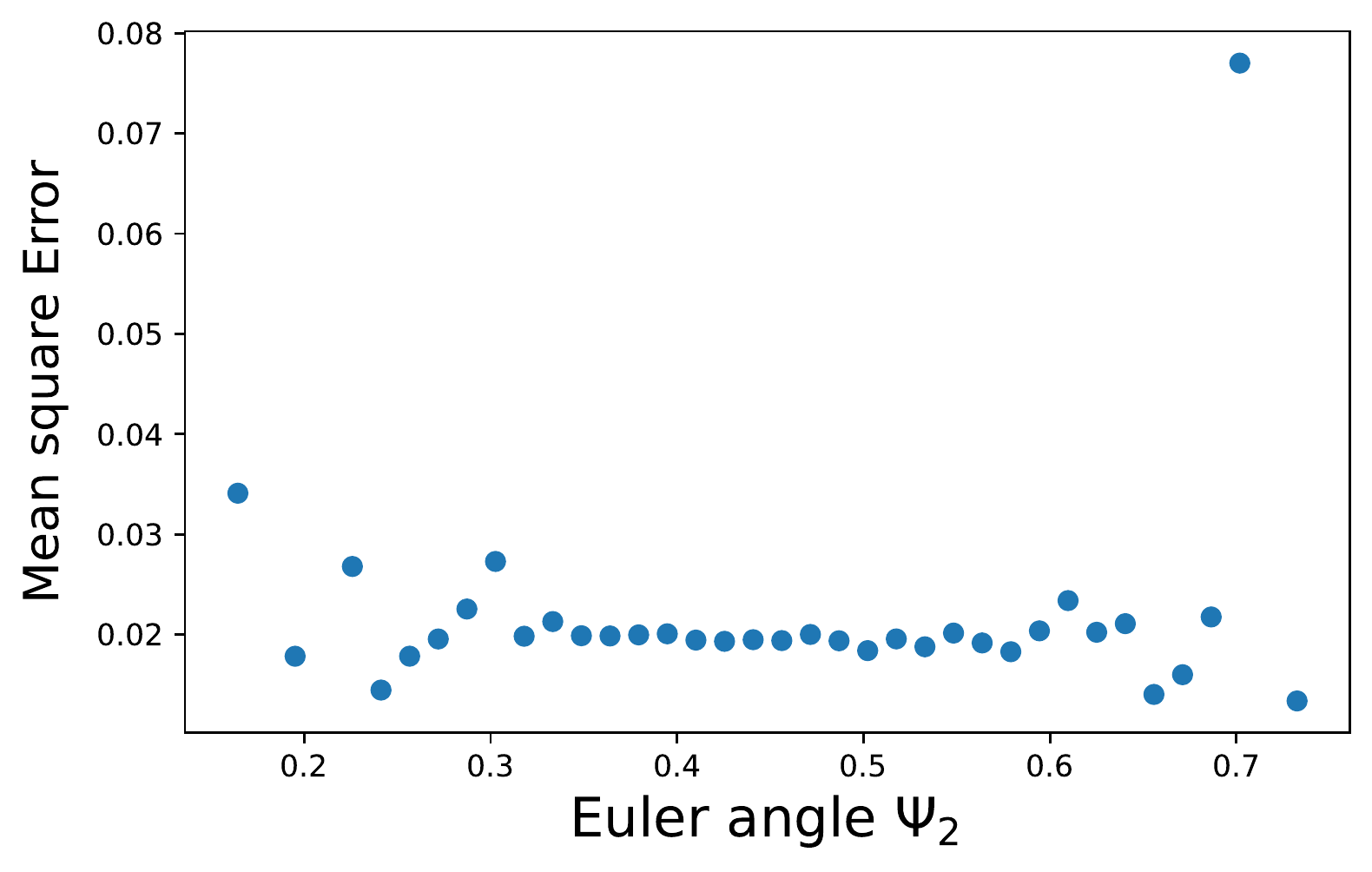}}
	\caption{Plot of average Euler angle against average MSE. 
		Each point represents a bin of Euler angle histogram. 
		The average of Euler angle and MSE is calculated over a bin. 
		Except for a few outliers, the variation of MSE with respect to Euler angle in all three directions is essentially constant. 
		The plot shows that model performance is independent of the Euler angle.}
	\label{fig:mse_euler}
\end{figure*}

Similarly, the effect of the Euler angles and the aspect ratio for each microstructure was compared against the MSE.
The effect of the aspect ratio on MSE is shown in Figure {\ref{fig:mse_aspect}}.
The prediction is better for microstructures with grains that are closer to circular rather than elongated.
Figure \ref{fig:mse_euler} shows that the variation of the MSE remains almost the same for Euler angles in all three directions except for small outliers. 
This suggests that Euler angles do not affect the model performance.

\subsection{Error Analysis of the  Predicted Stress Fields Using Cosine Similarity}

Since the MSE can range from $[0,\infty)$, it is difficult to use this to interpret the accuracy of the predicted stress fields when compared to the ground-truth.
Hence, the cosine similarity (Equation \ref{eq:cosine_simialrity}) is also used to analyze the model accuracy.
If $\mathbf{A}$ and $\mathbf{B}$ represent the predicted and the ground-truth stress field for a given microstructure, the cosine similarity is given by:
\begin{align}
	\text{Cosine similarity}(\mathbf{A},\mathbf{B}) &= \frac{\mathbf{A}\cdot \mathbf{B}}{\| \mathbf{A} \|\| \mathbf{B} \| } \label{eq:cosine_simialrity}
\end{align}

\begin{figure*}[htp]
	\centering
	\includegraphics[width=.45\textwidth]{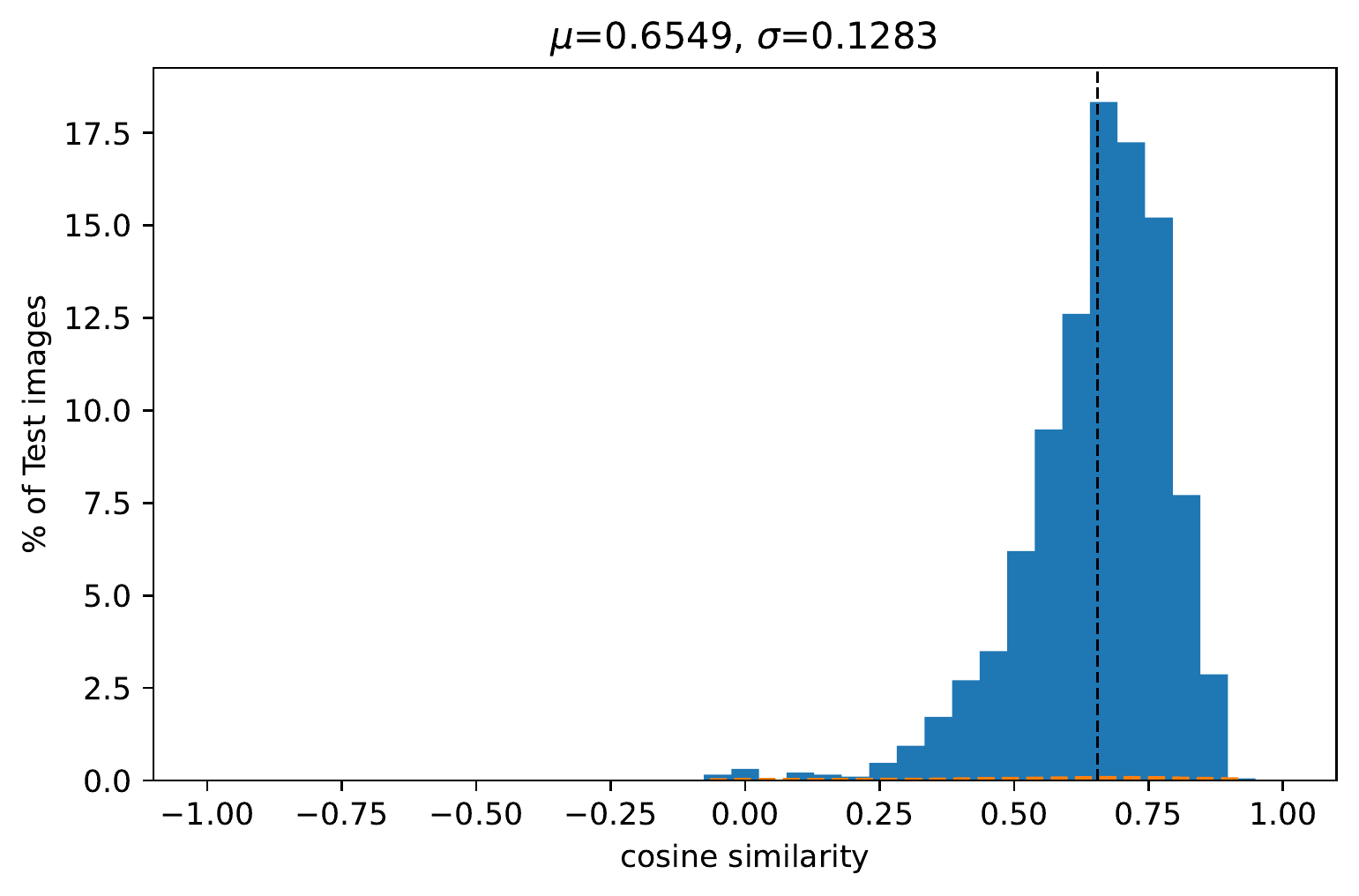}
	\caption{Histogram of cosine similarity measures for test images. 
		The horizontal and vertical axes represent the cosine similarity and percentage of samples respectively.
		The parameters $\mu$ and $\sigma$ represent the mean and the standard deviation respectively of the corresponding distributions. The cosine similarity ranges between $[-0.0908, 0.9217]$ indicating that for some microstructures, high-stress regions in predicted stress fields are overlapping with low-stress regions in the ground truth or vice versa.}
	\label{fig:dot_accuracy}
\end{figure*}

Before computing the cosine similarity, we scale the predicted and ground-truth stress fields between $[-1,1]$. 
The cosine similarity will range between $[-1,1]$, with the value of $1$ representing perfectly similar images. 
If the cosine similarity is computed when the pixel values of both the predicted stress field and the ground-truth range between $[0,1]$, as explained in the Section {\ref{sec:gen_preprocess}}, it would not provide sufficient information about the accuracy.  
If the low-stress regions in predicted stress fields coincide with peak-stress regions in the ground-truth, the cosine similarity will be low but is still positive.
By scaling between $[-1,1]$, we were able to distinguish the high-stress regions from the low-stress regions by making values of the low-stress regions negative.
If the high-stress regions in the ground-truth overlap with the low-stress regions in the predicted stress fields (or vice versa), the product will be negative, acting as a penalty to the cosine similarity. 
Further, the cosine similarity metric includes the norm of images in its denominator, which normalizes the stress values in each image. 
This helps in analyzing how accurately the hotspots are overlapping, irrespective of the hotspot magnitudes. 
In Figure {\ref{fig:dot_accuracy}}, it can be seen that some cosine similarity values are negative, which suggests that the predictions in such images are completely wrong. One such case is shown in the Figure {\ref{fig:worst_case}}. 
On further investigation, it was observed that the images with the negative cosine similarity values also had very high MSE. 
Hence this metric is a good indicator to compare the distribution of the stress fields irrespective of their magnitude.
The predicted stress fields may not be perfectly similar; however, high positive values indicate that peak-stress clusters are approximately overlapping, and for almost all images, peak-stress clusters do not overlap with the low-stress regions.

\subsection{Error Analysis of Clusters' Geometric Characteristics Inside Predicted Stress Fields}\label{sec:hotspots_analyse}

The accuracy of peak-stress cluster detection in the predicted stress fields is also studied by comparing their characteristics such as location, size, shape, and orientation with clusters in the ground-truth.  
For the current work, clusters around the three largest peaks in the microstructure are compared. 
As per our definition of clusters, the characteristics will change with the threshold, as shown in Figure \ref{fig:hotspot}. 
Threshold values of $\{0.5, 0.7, 0.8, 0.9\}$ were considered for our analysis. 
This analysis helps us understand what kind of clusters, depending on peak and threshold, can be detected inside the stress fields predicted by our Convolutional Encoder-Decoder model more accurately. 

\subsubsection{Comparing Location of Peak-Stress Clusters}

\begin{figure}[htp]
	\centering
	\includegraphics[width=.45\textwidth]{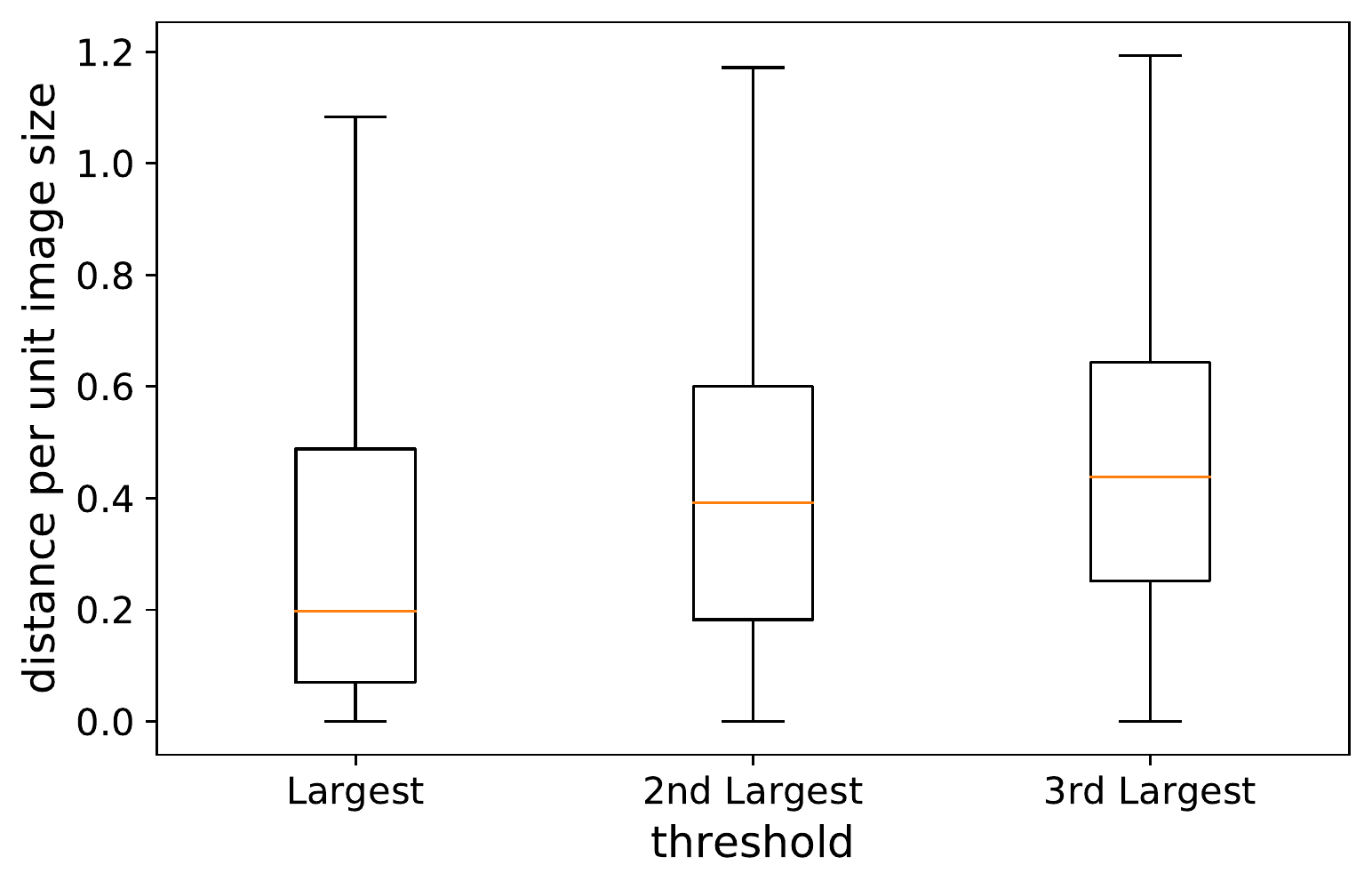}
	\caption{Each box plot represents the distribution of the distance (error in location) between the largest, second-largest, and third-largest peaks in the ground-truth and predicted stress fields.
		The model predicts the location of higher intensity peaks more accurately.
	}
	\label{peaks_dist}
\end{figure} 

The peak-stress cluster location is defined as the location of the peak-stress surrounded by the cluster, see Section \ref{sec:hotspots_characterize}; therefore, the cluster location is not affected by the threshold. 
The cluster location accuracy is analyzed using the Euclidean distance between stress peaks in the predicted stress fields and stress peaks of the same rank in the ground-truth. 
Figure \ref{peaks_dist} shows the distribution of distances in the test samples. 
It shows that the location of clusters around the largest peaks in predicted stress fields is more accurate. 
As the intensity of stress values in these clusters is large, the squared error for these peaks will be larger if the model predicts values inaccurately. 
Hence the model tries to learn by mainly focusing on the values at these peaks. 
Therefore, the prediction of the location of the highest peaks is much better.

\subsubsection{Comparing Size of Peak-Stress Clusters} 
\label{size-cluster}

\begin{figure}[htp]
	\centering
	\includegraphics[width=.5\textwidth]{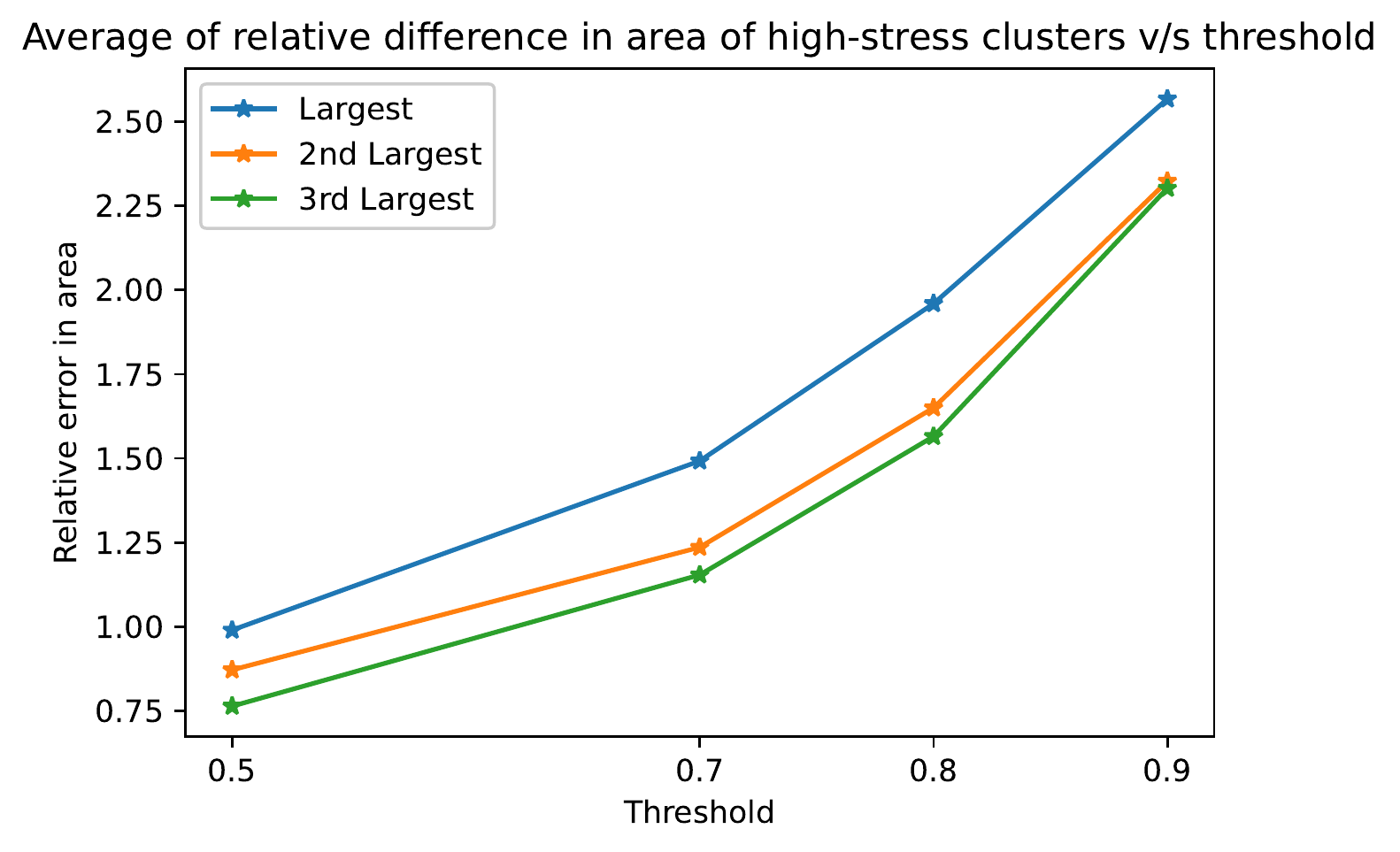}
	\caption{The mean absolute relative error of the areas of the detected clusters in the predicted stress fields for different thresholds. 
		The model predicts the size of larger clusters more accurately compared to smaller clusters.
		The errors are lower for clusters with smaller thresholds and around smaller peaks since their area is relatively larger in ground-truth stress fields.}
	\label{area_dist}
\end{figure}

The size of clusters in the predicted stress fields was compared with clusters in the ground-truth using the absolute relative difference in the area, $\delta_{a}$, as defined below. $\text{Area}_{\text{groundtruth}}$ refers to area of a peak-stress cluster in a ground-truth stress field. 
\begin{align}
	\delta_{a} = \left|\frac{ \text{Area}_{\text{groundtruth}}-\text{Area}_{\text{predicted}}}{\text{Area}_{\text{groundtruth}}}\right| \label{eq:ARDA}
\end{align}
As shown in Figure \ref{area_dist}, with an increase in the threshold, the relative difference increases.
Also, as the threshold increases, the clusters become smaller in size. 
Since, $\delta_{a}$ is inversely proportional to the size of the cluster in ground-truth, $\delta_{a}$ can be high for small clusters even with a small error in the detection of the number of pixels in the clusters in the predicted stress field as compared to the ground-truth.
As the higher threshold provides details about more concentrated clusters around the peak, the comparison shows that the model does not capture concentrated details around the peaks; instead, it captures the overall spread of the peak-stress clusters.

On comparing $\delta_{a}$ for clusters around the different peaks and same threshold, the model performs well where there are stress clusters around smaller peaks in the ground-truth. 
In the stress fields of most microstructures, these clusters are of moderate intensity but more spread out. 
Hence, the details of such clusters are predicted more accurately by the model. 
For similar reasons, the $\delta_{a}$ is lowest for the clusters around the third largest peaks and for all threshold values. 

The model's ability to predict the size of larger clusters more accurately can be attributed to the objective function $\MSE$.  
While training, the model emphasizes predicting a large number of pixels of stress clusters with average magnitude more accurately instead of a much smaller number of pixels but with relatively higher magnitude.

\subsubsection{Comparing the Orientation of Peak-Stress Clusters}

The accuracy of cluster orientation, $\theta$, was compared using the absolute difference, $\delta_{\theta} := \left|\theta_{\text{groundtruth}}-\theta_{\text{predicted}}\right|$, between the orientation of cluster in predicted stress fields and the ground-truth. 

Among the clusters defined using threshold values of $\{0.5,0.7,0.8\}$, the average $\delta_{\theta}$ is relatively small for the largest peak, as shown in Figure \ref{fig:theta}. 
This suggests that the model tries to capture the direction of the spread of clusters around higher peaks more accurately and gives lower priority to around lower peaks. 
This behavior can be attributed to the objective function MSE, where the model prioritizes learning high-intensity stress values.
For a threshold value of $\{0.9\}$, the behavior is the opposite; as the clusters are small, adding or removing even a single pixel in the clusters can change the orientation of the equivalent ellipse by a significant amount causing the high error. 

While this definition is very sensitive, this metric is still useful because it tells us that the orientation of clusters around larger peaks is more accurate and that the model is prioritizing clusters around the large peaks.
Due to the sensitivity of this metric, an unacceptable range would be $90$ degrees or more. 
Thus, for example, if a peak-stress cluster predicted by the model is oriented horizontally, whereas, in the ground truth, it is oriented vertically, this prediction would be completely unacceptable.
Other, less sensitive, metrics could provide a better picture of the orientation accuracy.

This analysis suggests that the orientation of clusters defined using a threshold of $0.7$ or less is fairly accurately predicted using our approach.

\subsubsection{Comparing the Shapes of Peak-Stress Clusters}

The shape was compared using the relative absolute difference in aspect ratios, $\delta_{ar} := \left|\frac{ \text{AR}_{\text{groundtruth}}-\text{AR}_{\text{predicted}}}{\text{AR}_{\text{groundtruth}}}\right|$, where $\text{AR}_{\text{groundtruth}}$ is the axis ratio of a peak-stress cluster in a ground-truth stress field.

\begin{figure}[htp]
	\centering
	\subfloat[comparison of orientation \label{fig:theta}]{\includegraphics[width=.45\textwidth]{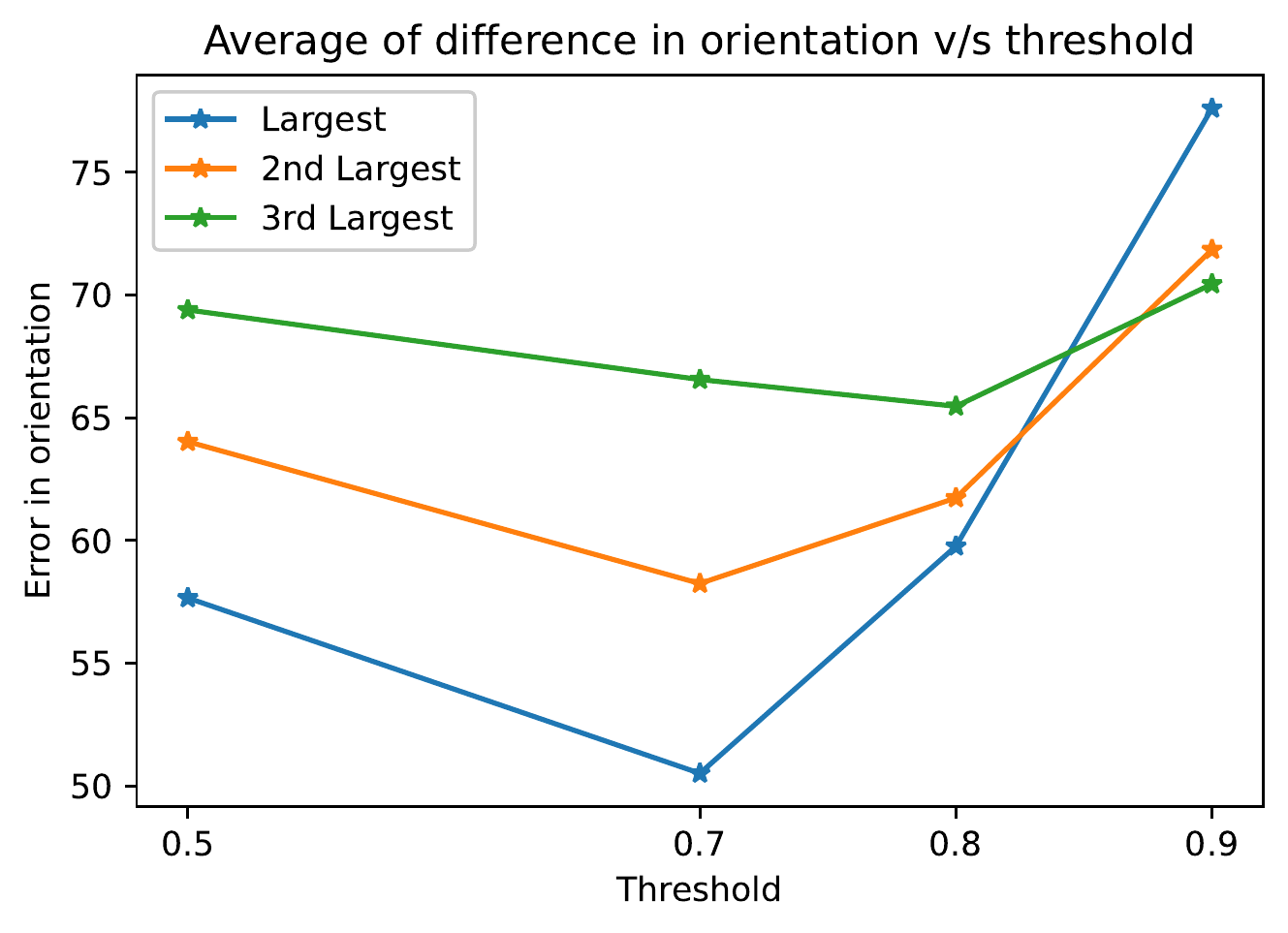}}
	\hfill
	\subfloat[comparison of aspect ratio \label{fig:axis}]{\includegraphics[width=.45\textwidth]{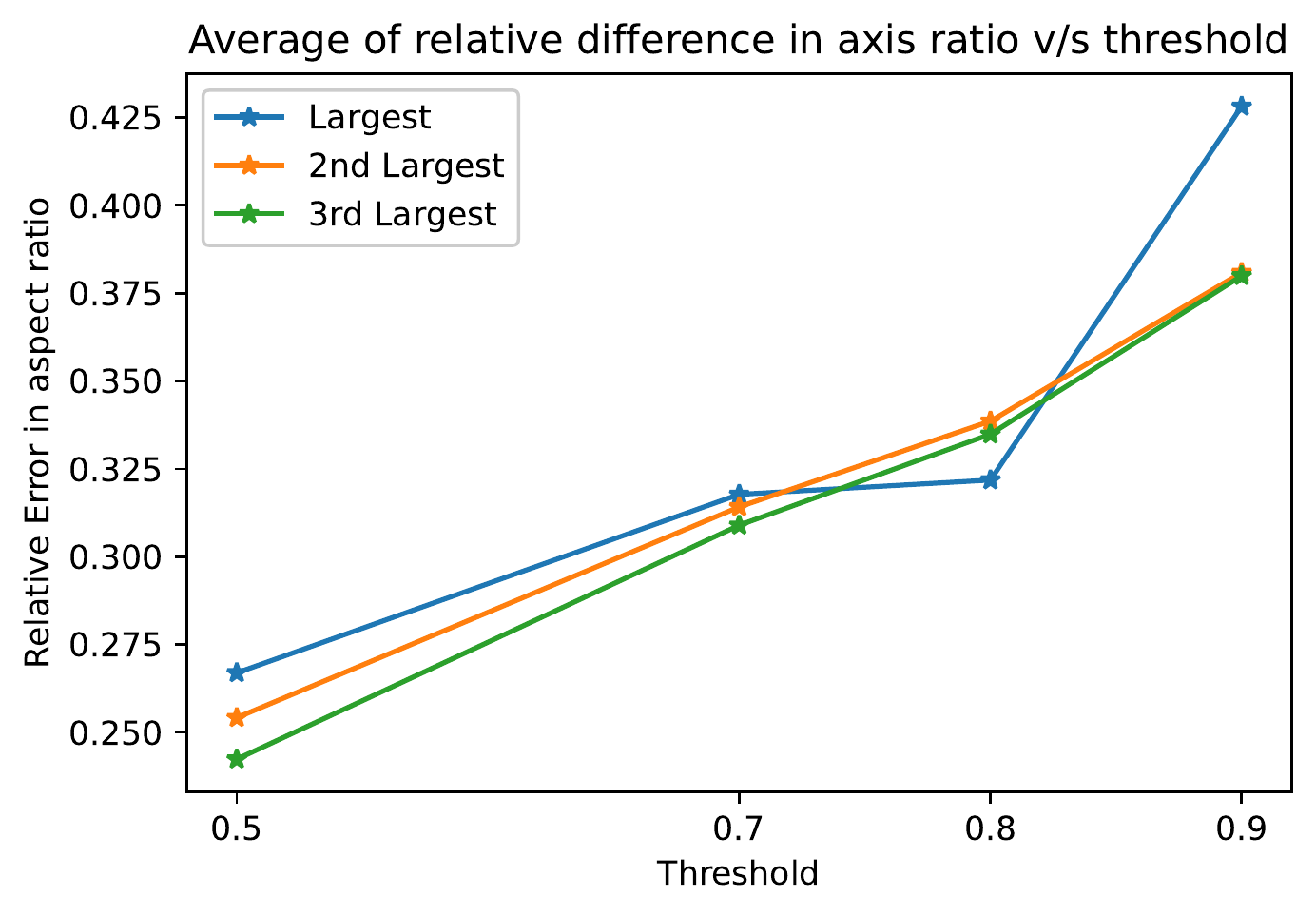}}
	\caption{Mean absolute difference of the orientation (a) and the mean absolute relative difference of the aspect ratio (b) of stress clusters in predicted and ground-truth stress fields for different thresholds. 
		The orientation and shape of clusters defined using the threshold of $0.7$ or less are fairly accurate, and clusters details around the largest peak are fairly accurate.}
\end{figure}

As shown in Figure \ref{fig:axis}, $\delta_{ar}$ is almost the same for all three clusters irrespective of the threshold used to define them.  
However, the error increases with an increase in the threshold. 
Since the clusters with lower thresholds are larger in size, for the same reason as in Section \ref{size-cluster}, the axis-ratio is more accurate for these clusters. 
With the increase in threshold, the size decreases, and even a slight shift of a pixel in clusters can change the aspect ratio by a considerable value hence the high error. 
This analysis suggests that the aspect ratio can be more accurately estimated for all three clusters defined using the threshold of $0.8$ or less using our approach. 

The complete analysis for cluster characteristics provides the understanding that the peak-stress clusters around the largest peak defined using a $0.7$ or less threshold can accurately be estimated using our approach. 
The stress regions around lower peaks and with a higher threshold may not be that accurate.

\subsection{Computational Cost of Convolutional Encoder Decoder Model}

The computational cost of the CED model was analyzed by estimating the time taken to train the model on the training set and the time taken by the trained model to predict the low-resolution stress fields for the complete microstructure dataset. 
A CPU machine with the configuration Intel(R) i7-3770 CPU @ 3.40 GHz was used to analyze the computational cost.

The machine took $8$ hours, $32$ minutes, and $14$ seconds to train the CED model using $8000$ training samples and training parameters as shown in Table {\ref{tab:optimisation}}. 
For predicting the low-resolution stress fields for $10000$ microstructures using the trained CED model, the machine took $1$ minute and $25$ seconds. 
Using the EVPFFT package {\cite{Lebensohn2012}}, the machine took $10$ hours, $34$ minutes, and $36$ seconds to simulate stress for the same $10000$ microstructures. 
The package was parallelized, and two processors were used to perform the simulations.

The CED model is around $448$ times more computationally efficient as compared to the EVPFFT package. 
However, these comparisons are limited only to CPU machines, and the efficiency can be different on GPU machines.
It is also important to note that training the machine learning model requires a relatively large memory compared to the EVPFFT package. 
This is because the model needs to load the complete dataset into RAM before training. 
A data-efficient training scheme can help reduce the requirement of large RAM.

\section{Feature Visualization}\label{Feature_visualisation}

The output of the filters from the first layer of the CED model was analyzed to investigate the microstructure features extracted by the model while predicting the stress fields.
In this work, the analysis was limited to just the first convolutional layer, as the filter behavior in the successive layers does not directly depend on the microstructure but is influenced by information extracted by the predecessor layers. 
As shown in the Table \ref{tab:arch}, there are $64$ filters in first layers. 
Each filter outputs a feature map for a given microstructure.

\begin{figure}
	\centering
	\includegraphics[width=.8\textwidth]{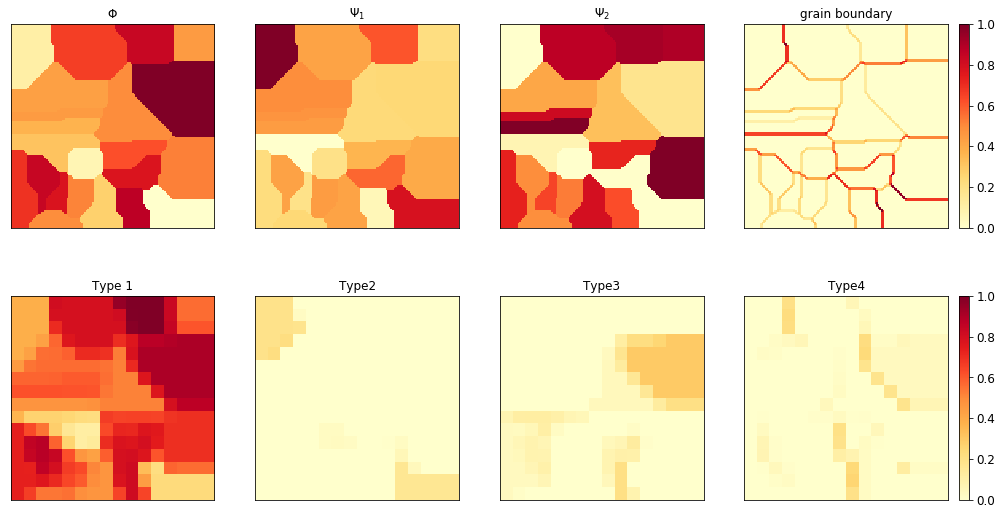}
	\caption{The first row shows the scaled input to the CED model. 
		The output of filters in the first layer of the model after activation was analyzed and grouped into five types based on their behavior. 
		The second row shows the output of selected filters belonging to the first four types. 
		The output from the filters belonging to the fifth type remains zero, and hence is not shown here. 
		Type 1 filters remain active in all regions. 
		Type 2 and Type 3 remain active only in regions with relatively high and low Euler angles, respectively.  
		Type 4 remains active only in regions where grain boundaries in specific orientation are present irrespective of Euler angle values.}
	\label{feature_map}
\end{figure} 

\begin{figure}[htb]
	\centering
	\includegraphics[width=.65\textwidth]{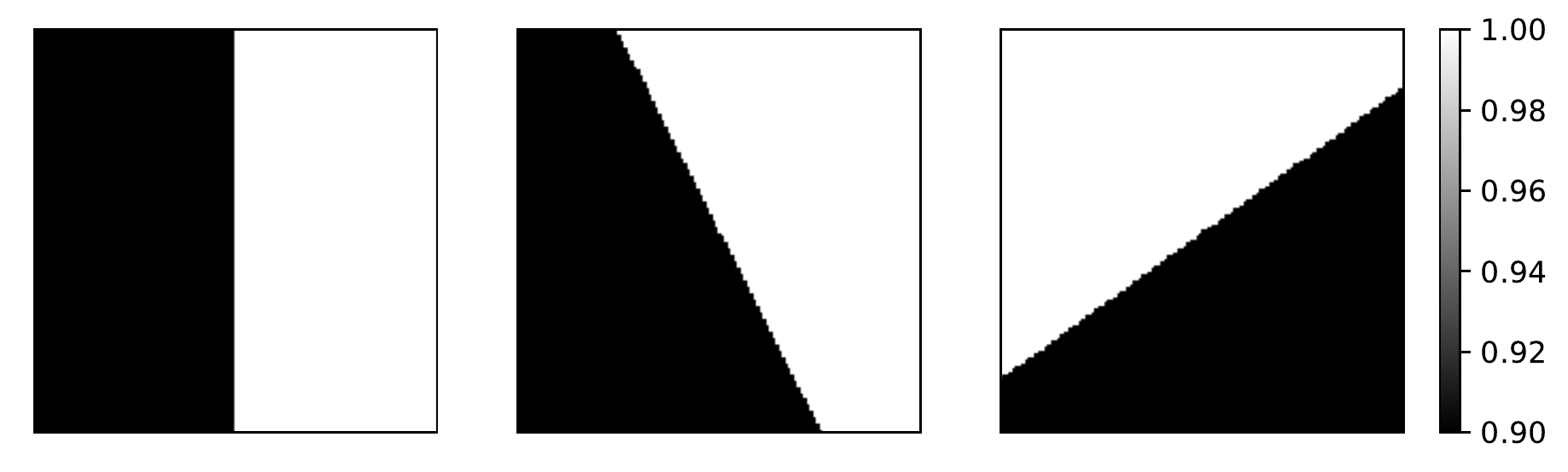}
	\caption{Manually designed bicrystal microstructures with the same values (normalized) of $(\Phi,\Psi_1,\Psi_2)$ but different grain boundary orientation. 
		The output of the filters in the first layer of the model was studied to analyze the role of grain boundary orientation.}
	\label{gb_rot}
\end{figure} 

Figure \ref{feature_map} shows the features extracted by selected filters of the model from a given microstructure and its grain boundary information. 
Since the microstructure had many complexities in such as grain boundaries, Euler angles, and grain junctions, we also study the behavior of these filters on simple manually-designed microstructures. 
In the first case, we considered a single crystal microstructure with the fixed Euler angle values. 
The values were varied from $0$ to $1$ since the scaled microstructure is fed into the model. 
In the second case, we considered a bicrystal with Euler angle values (normalized) but rotated the grain boundaries, as shown in Figure \ref{gb_rot}. 
The response of the filters was then studied. 
We classified filters based on their behavior as follows.

\begin{figure}[htb]
	\centering
	\includegraphics[width=.75\textwidth]{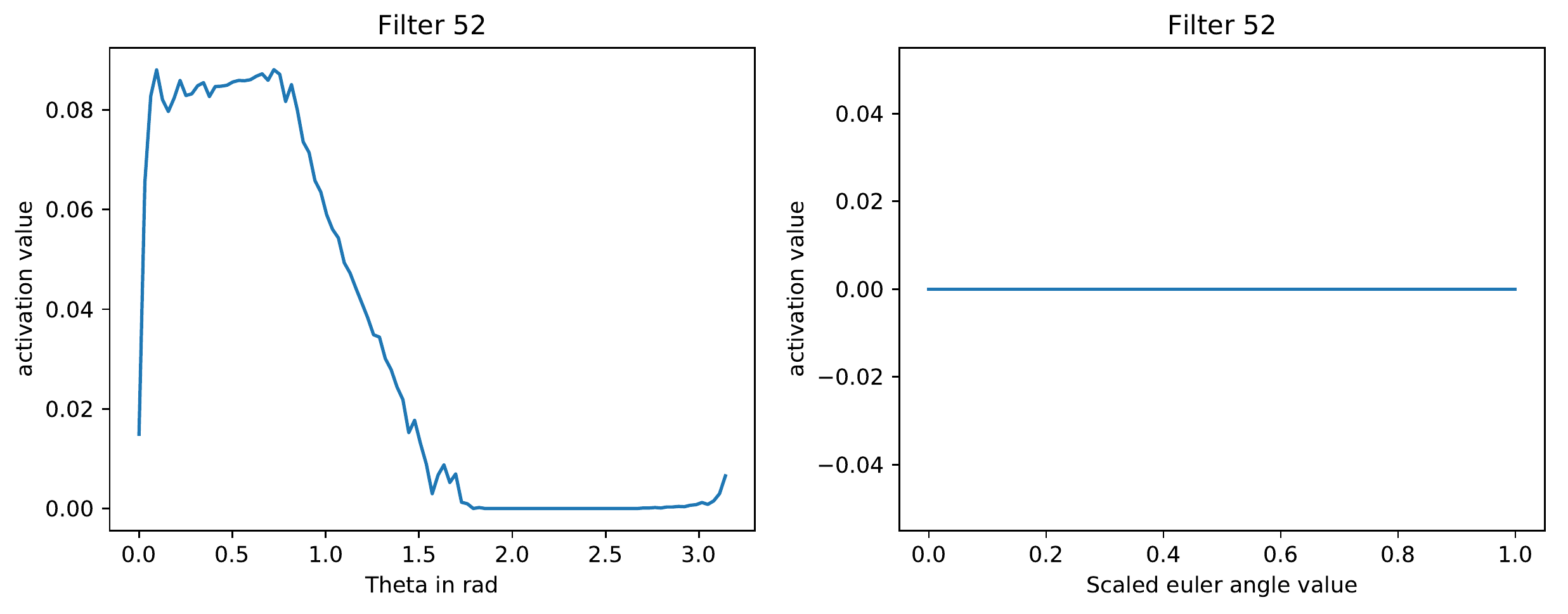}
	\caption{Activation plot of Filter 52 representing Type 4. The left figure shows activation of the filter when bicrystal microstructures with the same Euler angles but different orientations of grain boundary, as shown in Figure \ref{gb_rot} were fed to the model. On the right, the filter activation when monocrystalline microstructures with varying Euler angle values were provided to the model. These two figures show that the filter mainly focuses on detecting the grain boundaries and remain inactive in the regions with no grain boundary.}
	\label{type4}
\end{figure}

\begin{figure}[htb]
	\centering
	\includegraphics[width=.85\textwidth]{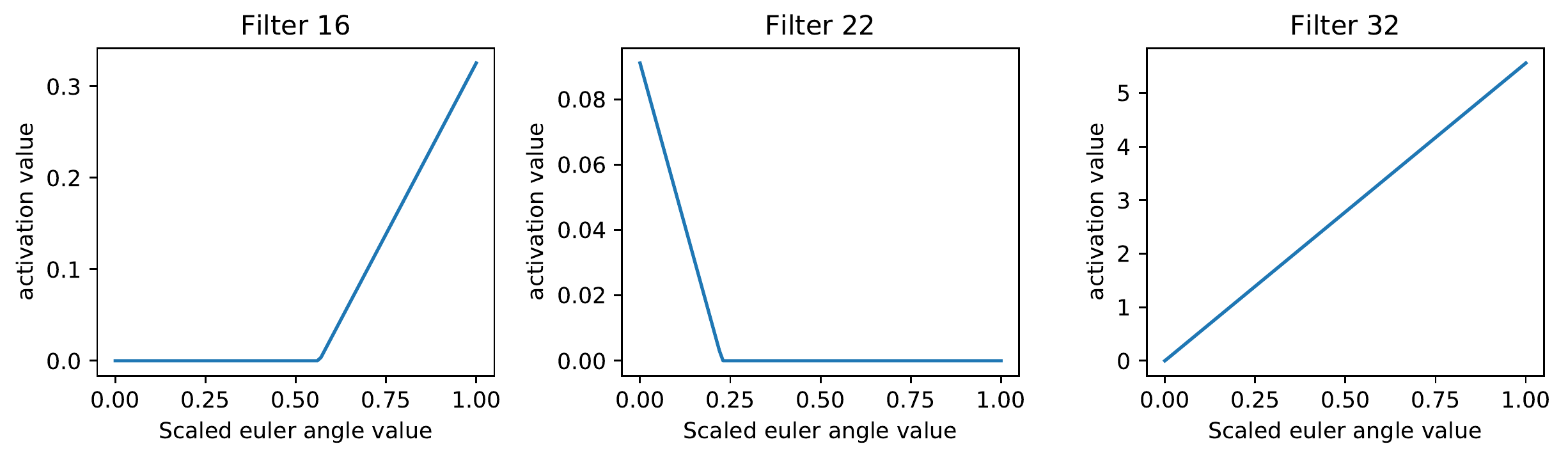}
	\caption{Activation plots with varying Euler angle values. Filter 16 represents Type 2, Filter 22 represents Type 3 and Filter 32 represents Type 1. Type 1 filters capture all the grains irrespective of their Euler angles, Type 2 filters detect grains with relatively larger Euler angles, and Type 3 filters detect grains with relatively smaller Euler angles.}
	\label{type_1_2_3}
\end{figure}

\begin{itemize}
	\item \textbf{Type 1}: These filters are active at every point in the microstructure irrespective of relative Euler angle values and grain boundary orientation. 
	The output after the $\relu$ activation from these filters will always be non-zero, as shown in Figure \ref{type_1_2_3}. 
	It means that these filters could extract information from all combinations of Euler angles and grain boundaries. 
	However, the magnitude of feature extraction varies with Euler angles and grain boundary orientation.
	
	\item \textbf{Type 2}: These filters are active only in regions with large Euler angles, Figure \ref{type_1_2_3}. 
	Their response decreases with a decrease in Euler angle values. 
	The response also varies with grain boundary orientation around grains with larger Euler angles. 
	They are unable to detect grain boundaries around grain with smaller Euler angles.
	
	\item \textbf{Type 3}: These filters are active only in the presence of small Euler angles, Figure \ref{type_1_2_3}. 
	Their response will decrease and then becomes $0$, with an increase in Euler angle values. 
	Their response varied with the orientation of grain boundaries but cannot detect grain boundaries around grains with large Euler angles. 
	It shows that the filters were able to extract features only with relatively small Euler angle values. 
	
	\item \textbf{Type 4}: These filters remain inactive for most cases but are active only in the presence of grain boundary in a specific orientation or due to significant relative difference between $\Phi,\Psi_1,$ and $\Psi_2$, as shown in Figure \ref{type4}. 
	It shows that the filters mainly focus on grain boundary detections, as shown in Figure \ref{feature_map}.
	
	\item \textbf{Type 5}: These filters remain inactive in every case.
	The parameters of these filters are negative; hence the output after the $\relu$ activation from these filters will always be zero.
\end{itemize}

\begin{figure}[htb]
	\centering
	\includegraphics[width=.5\textwidth]{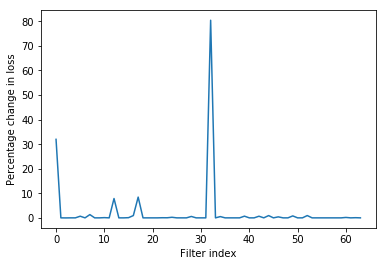}
	\caption{Percentage change in test error when a filter with a given index (plotted on the horizontal axis) was removed from the model. The plot shows that filter $0$ and $32$ play the most crucial role in stress field prediction as the increase in test error after removing these filters is higher than the other filters.}
	\label{act_error}
\end{figure} 

This observation is consistent with the observations in Figure \ref{feature_map}. 
Since the microstructure values before being fed to the model are scaled from $0$ to $1$, the behavior of filters explains that the model is concerned with only relative variations in microstructure features, and it does not care about the values. 
It is important to note that the size of the single-crystal and bi-crystal microstructures used to analyze the filter behavior is $8 \times 8$. It represents a small section of a large microstructure. 
As discussed previously, the model performance on single-crystals and bi-crystals of size $128 \times 128$ are highly inaccurate; hence visualization on such microstructures can lead to wrong interpretations. 
Moreover, the output of the filters in the first layer on the microstructure of the size $8 \times 8$ was a scalar value, which was easy to analyze.

Finally, we study the importance of each filter in predicting the stress distribution. 
For this, we select a filter, set its parameters to $0$, and predict the stress distribution on the test set. 
Setting the parameter to $0$ is equivalent to removing those filters, as the output from these filters will not influence the output of successive layers.
The test loss obtained after removing a particular filter from the model is compared against the original test loss. 
Figure \ref{act_error} shows the percentage change in loss. 
The error in the model prediction on test samples increased by almost $80\%$ on removing filter $32$.
This filter also remains active at every point in the microstructure, detecting every feature. 
Hence, this filter is the most important. Similarly, filter $0$ is the next most important filter. 

\section{Conclusion}\label{conclusion}

In this work, we presented a method to predict peak-stress clusters in polycrystalline microstructures under a macroscopic strain. 
First, stress fields were predicted using a Convolutional Encoder-Decoder model for given microstructures. 
A combination of a peak detection algorithm and a connected component labeling algorithm was used to detect peak-stress clusters in the predicted stress fields. 
The method was demonstrated for microstructures of pixel size $128 \times 128$, and peak-stress clusters were detected in the stress field of resolution $32 \times 32$ predicted by the CED model. 
The model was trained using $7000$ microstructure samples. 
The performance of the model was analyzed for $2000$ test samples by comparing predicted fields with the ground-truth using mean squared metrics for the overall stress field. The accuracy of detected peak-stress clusters was analyzed using the cosine similarity metric and comparing their geometric characteristics.

On investigating the mean squared error against the average grain sizes, average grain shapes, and average Euler angles of microstructures, it was observed that the CED model performed well for samples with average (normalized) grain size in the range of $[0.2,0.4]$.
The mean squared error for the microstructures with the average grains of size outside this range was high. 
The cosine similarity metric was used to analyze the overall accuracy of peak-stress clusters inside the predicted stress field. Further, the location, size, shape, and orientation of peak stress clusters of threshold values of $\{0.5,0.6,0.7,0.9\}$ was individually compared. 
This analysis showed that the geometric characteristics such as size and aspect ratio of the peak-stress clusters around the highest peaks are more accurate. The orientation of lower threshold clusters around larger peaks is more accurate, and the location of the clusters around the largest peak is most accurate.

The reason for this good performance is because the model prioritized learning the characteristics of very high-stress values. 
Further, peak-stress clusters defined using a threshold value $ \leq 0.7$ are more accurate than stress clusters with a threshold value of $0.7$ or more. 
The reason for this is that the model was able to better capture the large-scale structure of peak stress clusters and missed the finer-scale details.
The overall analysis provides insight into the microstructure features that are most important in determining the peak stresses.

The proposed CED model was trained and further evaluated to predict low-resolution von Mises stress fields. 
With a slight change in hyperparameters, the model can also be trained to predict other quantities such the principal stresses.
Although our method shows promising results, a CED model with a different architecture can potentially provide more accurate predictions.
For instance, a deep learning architecture that can predict stress fields with the same resolution as the microstructure, instead of a coarsened stress field, will likely be more effective in providing accurate predictions.
However, the dimensionality will be significantly larger, and hence training and prediction will be vastly more expensive.
Alternatively, a CED model with a fully-connected layer can be used to analyze the compressed information in the lower dimension \cite{guo2017deep}, to go beyond the current work wherein we did not use fully-connected layers as it suffered from early over-fitting.

An important issue that we did not address is the evolution of peak stresses under evolving load.
A CED model combined with Long Short Term Memory architecture can potentially be designed to predict peak-stress time evolution in large systems without significant recomputation at each time.

Finally, a very important challenge in machine learning-based models is not only to make predictions but also to understand the physics underpinning the predictions.
Various visualization algorithms such as the smooth-grad \cite{smilkov2017smoothgrad}, and guided back-propagation \cite{springenberg2014striving} can be used to investigate the specific features of the microstructure that play a significant role in the formation of peak-stress clusters.

\section*{Declaration of Conflicting Interests}

The authors declare that there is no conflict of interest.

\section*{Research Data}

The DREAM3D software to generate microstructure is available publicly at \url{dream3d.bluequartz.net/}. A version of the code developed for this work is available at  \url{github.com/ankitvaibhava/Peak_Stress_in_Microstructures}.

\section*{Funding}

The authors disclosed receipt of the following financial support for the research, authorship, and/or publication of this article: This work was supported by the National Science Foundation [DMREF 1921857 and CMMI MOMS 1635407], Army Research Office [MURI W911NF-19-1-0245], and Office Of Naval Research [N00014-18-1-2528].
 
\begin{acknowledgments}
    We thank Pittsburgh Supercomputing Center for XSEDE computing resources; 
    Ankita Mangal and Vahid Tari for help with learning to use Dream3D and EVPFFT;
    and Profs. Pei Zhang and Susu Xu, and Drs. Mostafa Mirshekari, Jonathon Fagert and Rajat Arora for useful discussions.
\end{acknowledgments}

\bibliographystyle{alpha}
\bibliography{ML-stress}

\end{document}